\input amstex.tex
\magnification=\magstep{1.5}
\baselineskip 22pt
\documentstyle{amsppt}
\topmatter
\title
Fundamental groups of open  $K3$  surfaces
Enriques surfaces and Fano 3-folds
\endtitle
\author
J. Keum and D. -Q. Zhang
\endauthor
\address
\endaddress
\abstract
We investigate when the fundamental group of the smooth part of a
$K3$  surface or Enriques surface with Du Val singularities,
is finite. As a corollary
we give an effective upper bound for the
order of the fundamental group of the smooth part of a certain Fano 3-fold.
This result supports Conjecture A below, while Conjecture A
(or alternatively the rational-connectedness conjecture in [KoMiMo]
which is still open when the dimension is at least 4)
would imply that every log terminal Fano
variety has a finite fundamental group.
\endabstract
\leftheadtext{J. Keum and D. -Q. Zhang}
\rightheadtext{Fundamental groups}
\endtopmatter
\document
\head
\endhead

\par \vskip 2pc
{\bf Introduction}

\par
We work over the complex numbers field  $\bold{C}$.
In this note, we consider the topological fundamental group
$\pi_1(T^0)$ of the smooth part $T^0$ of a normal projective
variety $T$. In general, it is difficult to calculate
such groups. Even in surface case, we still do not know
whether there is a plane curve $C$ such that the group
$\pi_1({\bold P}^2 \setminus C)$ is non-residually finite;
we note also that only in 1993, D. Toledo constructed
the first example of compact complex algebraic variety
with non-residually finite fundamental group,
which answered a question of J. P. Serre.

\par
In the present paper, the algebraic variety $T$ is assumed
to be either a $K3$  surface, or
an Enriques surface or a {\bf Q}-Fano 3-fold,
which has at worst log terminal singularities.
We will see from Theorem 3 and its proof that $\pi_1(T^0)$
of these three different objects are closely inter-related.

\par
First, let $X$ be a $K3$  surface with at worst Du Val singularities
(which is certainly log terminal; see [Ka]).
Then $X$ is still simply
connected (cf. [Ko1, Theorem 7.8]).
By [Ni1, Theorem 1], the number $c = \#(\text{\rm Sing} X)$
is bounded by 16, and if $c = 16$ then $\pi_1(X^0)$ is
infinite (cf. Remark 1.4). Recently, Barth [B1] has extended
this result in the following way : if each point in
$\text{\rm Sing} X$ is of Dynkin type $A_n$ ($n \ge 2$)
then $c \le 9$ and in the case $c = 9$,
$\pi_1(X^0)$ is infinite.
Our Theorem 1 below also implies that the condition $c = 16$
(resp. $c = 9$) in the result of Nikulin (resp. Barth)
is actually necessary and sufficient for $\pi_1(X^0)$
to be infinite (see Theorem 1 below for the precise statement).

\par
A similar result is obtained for the fundamental group
$\pi_1(W^0)$ of the smooth part of an Enriques surface
with at worst Du Val singularities (Theorem 2). In contrast with the
K3 case, $\pi_1(W^0)$ may not be abelian and may not
be $p$-elementary in the abelian case.

\par
One motivation behind this note is Theorem 3 below
in connection with the study of higher dimensional geometry
and an attempt to solve the conjectures below.
In what follows, a normal variety $V$ with at worst log
terminal singularities is {\bf Q}-Fano if, by definition, the anti-canonical
divisor $-K_V$ is {\bf Q}-Cartier and ample.

\par \vskip 0.8pc
{\bf Conjecture A.} {\it Let $V$ be a {\bf Q}-Fano $n$-fold.
Then the topological fundamental group $\pi_1(V^0)$ of
the smooth part $V^0$ of $V$ is finite.}

\par \vskip 0.8pc
{\bf Conjecture B.} {\it Let $V$ be a {\bf Q}-Fano $n$-fold.
Then the topological fundamental group $\pi_1(V)$ is finite.}

\par \vskip 0.8pc
{\bf Conjecture C.} {\it Let $V$ be a {\bf Q}-Fano $n$-fold.
Then $V$ is rational-connected.}

\par \vskip 0.8pc
Here $V$ is rational-connected, if any two general points
of $V$ can be connected by a single irreducible rational curve.
Clearly, Conjecture A implies Conjecture B.

\par
Conjecture A was proposed in [Z1] and was answered in affirmative
when the Fano index of $V$ is greater than $\dim V - 2$.
When $\dim V = 2$, Conjecture A was proved to be true
in [GZ1,2] or [Z2] (see [FKL] and [KM] for new proofs; see also [Z3]).

\par
Conjecture C implies Conjecture B [C, Ko1].
Conjecture C has been proved when $\dim V \le 3$ [C, KoMiMo],
but it is still open when $\dim V \ge 4$.
Our Theorem 3 below is a support towards Conjecture A.

\par \vskip 0.8pc
Now we state our Theorem 1.
Let $X$ be a $K3$  surface with Du Val singularities.
Let  $f : \widetilde{X} \rightarrow X$
be a minimal resolution,
$\Delta = f^{-1}(\text{\rm Sing}X)$
the reduced exceptional
divisor of  $f$  and  ${\bold Z}[\Delta]$  the sublattice
of  $H^2(\widetilde{X}, \bold{Z})$  generated by the cohomology classes
of irreducible components of  $\Delta$.
The universal covering map
$Y^* \rightarrow X^0 = X \setminus \text{\rm Sing} X$
can be extended to a morphism  ($Y^* \subseteq$)
$Y \overset{\gamma}\to{\rightarrow} X$
such that $Y/\pi_1(X^0) = X$; indeed,
if $\pi_1(X^0)$ is finite, then
$Y$  is the normalization of  $X$
in the function field  $\bold{C}(Y^*)$;
if  $\pi_1(X^0)$  is infinite, $\gamma$  is given
in Theorem 1 (3).

\par \vskip 0.8pc
{\bf Theorem 1.} {\it Let  $p$  be a prime number and
$X$ a  $K3$  surface with  $c$ ($c \ge 1$)
singularities of type  $A_{p-1} ($i.e., type $\frac{1}{p}(1, p-1))$
and no other singularities.
Then one of the $18$ rows in Table $1$ occurs; each of these
$18$ rows is realized by a concrete example.
Table $1$ shows precisely the topological
fundamental group $\pi_1(X^0)$  and
$\text{\rm Sing} Y$; in particular, we have:}

\par
(1) {\it $p \le 19$; if $p > 7$ then $\pi_1(X^0) = (1)$.}

\par
(2) {\it Suppose that
$\pi_1(X^0)$ is finite. Then  $\pi_1(X^0) = ({\bold Z}/(p))^k$
for some $0 \le k \le 4$ and  $Y$
$($a compactification of the universal cover of  $X^0)$
is a  $K3$  surface with at worst several type  $A_{p-1}$
singularities.}

\par
(3) {\it Suppose that $\pi_1(X^0)$  is infinite.
Then  $(p, c) = (2, 16)$ or
$(3, 9)$, and there is a  $\bold{Z}/(p)$-Galois cover
$X_1 \rightarrow X$ unramified over
$X^0$ such that  $X_1$ ($= \bold{C}^2/($a lattice$))$
is an abelian surface.
Hence we have an exact sequence:}
$$(1) \rightarrow \bold{Z} \oplus \bold{Z} = \pi_1(X_1)
\rightarrow \pi_1(X^0) \rightarrow \bold{Z}/(p) \rightarrow (1).$$
{\it The composition $\gamma$ of the natural morphisms
$Y = \bold{C}^2 \rightarrow X_1
\rightarrow X$, restricted over  $X^0$, is the
universal covering map of $X^0$.}

\par \vskip 0.8pc
Our next theorem utilizes Theorem 1 but needs some
lattice-theoretical arguments to determine the group
structure of the fundamental group.

\par \vskip 0.8pc
{\bf Theorem 2.} {\it Let  $p$  be a prime number.
Let $W$ be an Enriques surface containing a configuration
of smooth rational curves of Dynkin type $cA_{p-1}$
(the direct sum of $c \ge 1$ copies of $A_{p-1}$),
and let $W^0$ be the surface with these $c(p-1)$ curves on $W$ removed.

\par
Then one of the $26$ rows in Table $2$ occurs; in particular,
$\pi_1(W^0)$ is soluble and it is infinite if and only if $(p, c)
= (2, 8)$. Though the realization of the $2$ rows in Table $2$ are
unknown yet, each of the remaining $24$ rows in Table $2$ is
realized by a concrete example.}

\par \vskip 0.8pc
The following is an application of Theorems 1 and 2 and a
partial answer to Conjecture A above.

\par \vskip 0.8pc
{\bf Theorem 3.} {\it Let $p$ be a prime number.
Let $V$ be a Fano $3$-fold with a Cartier divisor
$H$ such that $m(K_V + H)$ is linearly equivalent
to zero for $m = 1$ or $m = 2$. Suppose that a member
$H$ of $|H|$ is irreducible normal and has $c$ singularities of
type $A_{p-1}$ and no other singularities.

\par
Then the fundamental group $\pi_1(V^0)$ of the smooth part $V^0$
of $V$ is the image of a group in Table $1$ or $2$.
In particular, $\pi_1(V^0)$ is soluble; and if
$(p, c) \ne (2, 8), (2, 16), (3, 9)$,
then $|\pi_1(V^0)| \le 2p^k$ for some $0 \le k \le 4$.}

\par \vskip 0.8pc
{\bf Remark 4.} (1) On a {\bf Q}-Fano
3-fold $V$, a relation $m(K_V + H) \sim 0$
with $H$ a Cartier divisor occurs when $V$
has Fano index 1 and Cartier index $m$.
It is conjectured that in this situation
$m = 1, 2$. This conjecture
is confirmed by T. Sano [Sa] under the stronger condition
that $V$ has at worst terminal cyclic quotient singularities.
On the other hand, a result of Minagawa [Mi] shows that
any terminal {\bf Q}-Fano 3-fold of Fano index 1 can be deformed to
a {\bf Q}-Fano 3-fold of Cartier index 1, 2.

\par
(2) By Ambro [A, Main Theorem], a general member of $|H|$
is normal irreducible and has at worst log
terminal singularities; so $H$ has at worst Du Val or
type $(-4)$ or type $(-3)-(-2)-\cdots-(-2)-(-3)$
singularities since $2K_H \sim 0$ (cf. the proof of Theorem 3), whence
the condition on Sing $H$
in Theorem 3 is quite reasonable.
By the proof in \S 4, we always have a
surjective homomorphism $\pi_1(H^0) \rightarrow \pi_1(V^0)$,
where $H^0 = H -$ Sing $H$.
In [SZ], a sufficient condition for $\pi_1(H^0)$
to be finite is given when $K_H \sim 0$.

\par
(3) The author has not been able to
construct an example of $V$ in Theorem 3 satisfying
$(p, c) = (2, 8), (2, 16)$ or $(3, 9)$.

\par \vskip 0.8pc
Theorems 1, 2 and 3 are proved respectively in \S 2,
\S 3 and \S 4.

\par \vskip 0.8pc
{\bf Acknowledgment.} Part of this work was done when
the second author was visiting University of Missouri and
University of Michigan in 1999. He would like to
thank both institutions for their hospitality.
This work was finalized during the first author's
visit to Singapore under the Sprint programme of
the Department of Mathematics. He thanks NUS
for the financial support and the hospitality.
Both authors are very thankful to Professor I. Shimada for the interest
in this paper who actually found an alternative
shorter proof of the Table 1 using computer.
They also thank Professor R. V. Gurjar for his crucial help
in proving Theorem 3. This work was partially financed
by an Academic Research Fund of National University of Singapore.

\par \vskip 2pc
{\bf \S 1. The  $K3$  case with $p = 2$}

\par \vskip 0.8pc
{\bf 1.1.} We will frequently and implicitly use the
following observation:
Let $f : \widetilde{X} \rightarrow X, \, c$ be as in Theorem 1.
Then there is a Galois $\bold{Z}/(p)$-cover
$\sigma : Z \rightarrow X$ ramified exactly over
a $c_1$-point subset $H = \{p_1, \dots, p_{c_1}\}$
of $\text{\rm Sing} X$
if and only if there is a relation
$\sum_{i=1}^{c_1} \Delta_i^* \sim p L$ on $\widetilde{X}$,
where $L$ is a Cartier divisor and $\Delta_i^*$ is an effective
Cartier divisor with support
equal to $\Delta_i := \cup_{i=1}^{c_1} f^{-1}(p_i)$ and coefficients
in $\Delta_i^*$ coprime to $p$.

\par
Indeed, assuming the above equivalent conditions,
one note that $\Delta_i = \sum_{k=1}^{p-1} \Delta_i(k)$
is a linear chain of $(-2)$-curves and can check that
$\Delta_i^* = d_i \sum_{k=1}^{p-1} k \Delta_i(k)$
for some integer $d_i$ coprime to $p$. Moreover,
one has $\Cal{O}(\overline{L})^{\otimes p} \cong {\Cal O}_X$
with $\overline{L}$ the $f$-image of $L$, and
$\sigma$ is given by:
$$\sigma : Z = Spec \oplus_{i=0}^{p-1} {\Cal O}_X(-i \, \overline{L})
\rightarrow X.$$

\par \vskip 0.8pc
{\bf Lemma.} (1) {\it Let $\Sigma \subset \text{\rm Sing} X$ and let
$X_1 \rightarrow X$ be a minimal resolution of singular
points not in $\Sigma$. Then $\pi_1(X_1^0) = \pi_1(X \setminus \Sigma)$,
where $X_1^0 := X_1 \setminus \text{\rm Sing} X_1$.}

(2) {\it One has
$H_1(X \setminus H, \bold{Z}) \cong
\overline{\bold{Z}[\sum_i \Delta_i]}/\bold{Z}[\sum_i \Delta_i]
\cong (\bold{Z}/(p))^{\oplus k}$ for some $0 \le k \le c_1$,
where for a sublattice $\Gamma$ of $H^2(X, \bold{Z})$,
we denote by $\overline{\Gamma}$ its primitive closure.}

\par
(3) {\it $\bold{Z}[\sum_i \Delta_i]$ is primitive in $H^2(X, \bold{Z})$
$\Leftrightarrow$ $\pi_1(X \setminus H)$
is a perfect group $\Leftrightarrow$ $H$ does not include any
$p$-divisible subset (cf. $1.2$ below).}

\par \vskip 0.8pc
{\it Proof.} (1) follows from [Ko1, Theorem 7.8] since $X$ has
at worst log terminal singularities.
The first isomorphism in (2) follows from the proof of
[X, Lemma 2], while the second follows from the
assumption on $\text{\rm Sing} X$. (3) is a consequence
of (1) and 1.1.

\par \vskip 0.8pc
{\bf Definition and Remark 1.2.} Let ${X}$ be as in Theorem $1$.
A subset $H$ of $\text{\rm Sing} {X}$ is $p$-{\it divisible}
if there is a Galois $\bold{Z}/(p)$-cover
${Z} \rightarrow {X}$ ramified exactly
over $H$ (cf. 1.1). When $p = 2$, $H$ is $2$-divisible
if and only if $f^{-1}(H)$ is $2$-divisible in
the lattice $Pic {\widetilde X}$ (see also 3.8).

\par \vskip 0.8pc
{\bf Lemma 1.3.} (cf. [Ni1, Lemma 3])
{\it Let $p, X, X^0, c$ be as in Theorem $1$.
Suppose that there is a Galois $\bold{Z}/(p)$-cover
$\sigma : Z \rightarrow X$,
ramified exactly over
$\text{\rm Sing} X$
$($i.e., $\text{\rm Sing} X$ is $p$-divisible$)$.
Then  $(p, c)$  fits one of the following cases:}
$$(2, 8), (2, 16), (3, 6), (3, 9), (5, 4), (7, 3).$$
{\it Moreover, if  $(p, c) = (2, 16), (3, 9)$,
then  $Z$  is an abelian surface and hence
$$\pi_1(X^0)/(\bold{Z} \oplus \bold{Z}) = \bold{Z}/(p);$$
if  $(p, c)$  fits one of the remaining $4$ cases,
then  $Z$  is a (smooth)  $K3$  surface
and hence $\pi_1(X^0) = \bold{Z}/(p)$.}

\par \vskip 0.8pc
{\it Proof.} By the assumption, for
each singular point $p_i$ of $X$, $q_i = \sigma^{-1}(p_i)$
is a single point and  $Z$  is smooth.
Now  $K_{X} \sim 0$  implies that
$K_{Z} \sim 0$, whence  $Z$  is either
abelian with Euler number  $e(Z) = 0$
or $K3$  with $e(Z) = 24$.  The lemma
follows from the calculation (noting that
K3 surfaces are simply connected):
$$e(Z) - c = p e(X^0)
= p(24 - c p).$$

\par \vskip 0.8pc
{\bf Remark 1.4.} There is a converse to Lemma 1.3
by [Ni1] and [B1].  Suppose that ${X}$
is a $K3$  surface with $\text{\rm Sing} {X}
= c A_1$ where $c \ge 16$ (resp. $\text{\rm Sing} {X}
= c A_2$ where $c \ge 9$). Then $c = 16$ (resp. $c = 9$)
and Sing $X$ is $p$-divisible with $p = 2$ (resp. $p = 3$);
so there is a Galois $\bold{Z}/(p)$-cover
${Y} \rightarrow {X}$  unramified over $X^0$
so that ${Y}$ is
an abelian surface. In particular, $\pi_1(X^0)$
is infinite soluble and all assertions in Theorem 1 (3) hold.
When $p = 2$, the covering involution
of ${Y}$ coincides with $\iota : (x,y) \mapsto (-x,-y)$.

\par \vskip 0.8pc
{\bf Lemma 1.5.} {\it Let $f : \widetilde{X} \rightarrow {X}, \, c$
be as in Theorem $1$ with $p = 2$.
Suppose that $H_1, H_2$ are two distinct $2$-divisible
$8$-point subsets of $\text{\rm Sing} {X}$.
Then either $c = 16$ and $\text{\rm Sing} {X} = H_1 \cup H_2$,
or $c \ge 12$ and $|H_1 \cap H_2| = 4$.}

\par \vskip 0.8pc
{\it Proof.} By 1.1, one has $D_i/2 \in H^2(\widetilde{X}, \bold{Z})$,
where $D_i = \sum_{x \in H_i} f^{-1}(x)$. Then
$(D_1+D_2)/2 \in H^2(\widetilde{X}, \bold{Z})$.
Now the lemma follows from 1.1 and Lemma 1.3 (noting that $c \le 16$
always holds by Remark 1.4).

\par \vskip 0.8pc
{\bf Lemma 1.6.} {\it Let ${X}, c$ be as in Theorem $1$
with $p = 2$.}

\par
(1) {\it Suppose that $c = 13$. Then there are $2$-divisible
$8$-point subsets $H_1, H_2$ of $\text{\rm Sing} {X}$
with $|H_1 \cap H_2| = 4$.}

\par
(2) {\it Suppose that $c \ge 14$. Then there is a $12$-point subset
$\Sigma$ of $\text{\rm Sing} {X}$ such that
$\Sigma$ includes only one $2$-divisible subset $H$.}

\par
(3) {\it Suppose that $c \ge 14$. Then there is an $11$-point subset
$\Sigma_1$ of $\text{\rm Sing} {X}$ such that
$\Sigma_1$ does not include any $2$-divisible subset.}

\par \vskip 0.8pc
{\it Proof.} In view of Lemma 1.5, for (1), it suffices to show
that $\text{\rm Sing} {X}$ includes two distinct
$2$-divisible (8-point) subsets. By the proof of [Ni1, Lemma 4]
or [B1, Lemma 2], $\text{\rm Sing} {X}$ includes
a $2$-divisible 8-point subset $H_1$. The same reasoning shows
that any 12-point set consisting of 7 points in $H_1$ and
the $5$ singular points of ${X}$ not in $H_1$,
includes a $2$-divisible 8-point subset $H_2$ ($\ne H_1$).
(1) is proved.

\par
For (2), applying (1), we get $2$-divisible
8-point subsets $H_1, H_2$ of $\text{\rm Sing} {X}$
with $|H_1 \cap H_2| = 4$. Take two singular points $p_1, p_2$
of ${X}$ not in $H_1 \cup H_2$, and one point
$p_3$ in $H_1$ but not in $H_2$. Applying Lemma 1.5, we see
that we can take $(H_1 - \{p_3\}) \cup H_2 \cup \{p_i\}$
as $\Sigma$, for $i = 1$ or $i = 2$.

\par
For (3), we let $\Sigma_1$ be any subset of $\Sigma$ in
(2) containing not more than 7 points of $H$.

\par \vskip 0.8pc
{\bf 1.7.} Let  $A = \bold{C}^2/\Lambda_A$  be an abelian surface
with $\iota$ the involution  $(x,y) \mapsto (-x,-y)$.
Denote by  $A_2$ the set of the 16 $\iota$-fixed points,
which is a subgroup of  $A$
consisting of the 2-torsion points. One can regard
$A_2$ as a $4$-dimensional vector space over
the field $\bold{Z}/(2)$.
The quotient ${X} := A/\langle \iota \rangle$
is a  $K3$  surface with 16 singularities $p_i$
of Dynkin type $A_1$ dominated by the points in $A_2$.
The bijection $A_2 \rightarrow \text{\rm Sing} {X}$
defines on the latter a $4$-dimensional $\bold{Z}/(2)$-vector
space structure.
One sees easily that $\pi_1(X^0)/\pi_1(A) = \bold{Z}/(2)$
and $\pi_1(X^0)$ is generated by the involution $\iota$
and  $\Lambda_A$.

\par
Suppose that $H$ is a $2$-divisible 8-point subset of
$\text{\rm Sing} {X}$ and
$\sigma : {Z} \rightarrow {X}$
the corresponding $\bold{Z}/(2)$-cover ramified
exactly over $H$. Then each singular point of
${X}$ not in $H$ splits into two type $A_1$
singularities of ${Z}$ and these 16 points
form the singular locus $\text{\rm Sing} {Z}$.
So ${Z} = B/\langle \iota \rangle$
with $B = \bold{C}^2/\Lambda_B$ an abelian surface
(Remark 1.4), and $\pi_1(Z^0)$ is generated by the involution
$\iota$ and  $\Lambda_B$.

\par
The covering $\sigma$ induces $\pi_1(X^0)/\pi_1(Z^0) = \bold{Z}/(2)$.
One can verify that $\Lambda_B$ is an index-2 sublattice of $\Lambda_A$.
This way, we obtain a commutative diagram:
$$\gather \bold{C}^2/\Lambda_B = B \longrightarrow
Z = B/\langle \iota \rangle \\
\hat{\sigma} \downarrow \hskip 2.2pc \downarrow \sigma \\
\bold{C}^2/\Lambda_A = A \longrightarrow
X = A/\langle \iota \rangle
\endgather $$
Note that $\hat{\sigma} : B_2 \rightarrow A_2$
(and hence $\sigma : \text{\rm Sing} Z
\rightarrow \text{\rm Sing} X$)
is a rank-3 linear map between $\bold{Z}/(2)$-vector
spaces of dimension 4.

\par \vskip 0.8pc
{\bf Lemma 1.8.} {\it Let ${X} = A/\langle \iota \rangle$
be a Kummer surface. Then we have:}

\par
(1) {\it An $8$-point subset $H$ of
$\text{\rm Sing} X$ is $2$-divisible if and only if
it is an affine hyperplane of the $\bold{Z}/(2)$-vector
space $\text{\rm Sing} X$.}

\par
(2) {\it For both $i = 1, 2$, there is a $12$-point subset $\Sigma_i$
of $\text{\rm Sing} X$ such that
$\pi_1(X_i^0) = \pi_1(X \setminus \Sigma_i)$ is equal to
$\bold{Z}/(2)$ (resp. $(\bold{Z}/(2))^{\oplus 2}$)
when $i = 1$ (resp. $i = 2$), where $X_i \rightarrow X$
is a minimal resolution of singularities not in $\Sigma_i$
and $X_i^0$ is the smooth part of $X_i$.}

\par \vskip 0.8pc
{\it Proof.} (1) follows from [Ni1, Cor. 5 and Remark 1].

\par
(2) Let $H_1, H_2$ be $2$-divisible 8-point subsets
of $\text{\rm Sing} X$ with $|H_1 \cap H_2| = 4$
(Lemma 1.6). As in 1.7, let
$\sigma : Z \rightarrow X$
be the double cover ramified exactly over $H_1$. Then
one can verify that $\hat{H}_2 := \sigma^{-1}(H_2) \cap \text{\rm Sing} Z$
$= \sigma^{-1}(H_2 \setminus H_1)$ is an affine hyperplane of
$\text{\rm Sing} Z$ and hence the group
$\pi_1(Z \setminus \hat{H}_2)$
equals $\bold{Z}/(2)$ (see the proof of Lemma 1.3).
The covering map $\sigma$ implies that
this group is an index-2 subgroup of
$\pi_1(X \setminus \Sigma_2)$ where $\Sigma_2 = H_1 \cup H_2$.
Hence the group $\pi_1(X \setminus \Sigma_2)$ has order 4;
since $X$ has no type $A_3$ singularity, this group equals
$(\bold{Z}/(2))^{\oplus 2}$ (cf. [X, Theorem 3]).
One has $\pi_1(X_2^0) = \pi_1(X \setminus \Sigma_2)$ by Lemma 1.1.

\par
By Lemma 1.6, we can find a 12-point subset $\Sigma_1$
of $\text{\rm Sing} X$ so that
$\Sigma_1$ contains only one $2$-divisible 8-point
subset $H$. As in 1.7, let
$\sigma : Z \rightarrow X$
be the double cover ramified exactly over $H$.
Let $g : \widetilde{Z} \rightarrow Z$ be a minimal resolution
with $\Gamma = g^{-1}(\text{\rm Sing} Z)$,
a disjoint union of 16 smooth rational curves.
The covering map $\sigma$ implies that $\pi_1(X \setminus \Sigma_1)$
has the (trivial) group in (iii) below as an index-2 subgroup.
So (2) is reduced to the proof of the claim below.

\par \vskip 0.8pc
{\bf Claim 1.8.1.} (i) The 8-point set
$\sigma^{-1}(\Sigma_1 \setminus H)
= \sigma^{-1}(\Sigma_1) \cap \text{\rm Sing} Z$
is not an affine hyperplane of $\text{\rm Sing} Z$,
and hence does not include any $2$-divisible set.

\par
(ii) The fundamental group of $Z$ with the
8 points in (1) removed, is trivial.

\par \vskip 0.8pc
If the first assertion of the claim is false, then
the 8-point set would be an affine hyperplane
and hence its $\sigma$-image is contained in an affine hyperplane
$H_3$ of $\text{\rm Sing} X$ which has to
consist of the 4 points $\Sigma_1 \setminus H$
and 4 points in $H$, and $\Sigma_1$ would include two
distinct $2$-divisible subsets $H, H_23$, a contradiction.

\par
By (i) and Lemma 1.1, the group in (ii) is perfect.
Moreover, this group is soluble and hence trivial because
it is the image of $\pi_1(Z^0)$ while
the latter is soluble [Remark 1.4]. This proves the claim
and also the lemma.

\par \vskip 0.8pc
{\bf Lemma 1.9.} {\it Let $f : {\widetilde X} \rightarrow X, \,
\Delta, \, c$ be as in Theorem $1$
with $p = 2$. Suppose that $\bold{Z}[\Delta]$ is primitive in
$H^2(\widetilde{X}, \bold{Z})$
(this is true if $c \le 7$; see Lemmas $1.1$ and $1.3$).
Then $c \le 11$ and $\pi_1(X^0) = (1)$ (each $c \le 11$ is realizable).}

\par \vskip 0.8pc
{\it Proof.} By the proof of [Ni1, Lemma 4] or [B1, Lemma 3],
the primitivity of $\bold{Z}[\Delta]$ implies that $c \le 11$.
>From [Ni2, Theorem 1.14.4] and its remark, one deduces that there is
a unique primitive embedding of $\bold{Z}[\Delta]$ into the
K3 lattice. Now by the connectivity theorem [Ni3, Theorem 2.10],
we are reduced to show the lemma for any particular ${X}_1$
satisfying the same condition of the lemma.

\par
Let $\widetilde{X}$ be a (smooth) Kummer surface with 16 disjoint
smooth rational curves. By Lemmas 1.1 and 1.6, among these 16,
there are 11 curves $E_i$ ($1 \le i \le 11$) such that
if $\widetilde{X} \rightarrow {X}_1$ is the contraction of
$E_1, \dots, E_c$ ($c \le 11$) then $\text{\rm Sing} X$ does not include
any $2$-divisible subsets.
Thus $\pi_1(X_1^0) = (1)$ as in the proof of Claim 1.8.1.
The lemma is proved.

\par \vskip 0.8pc
{\bf Proposition 1.10.} {\it Let
$f : \widetilde{X} \rightarrow X$,
$\Delta$, $c$ be as in Theorem $1$
with $p = 2$.}

\par
(1) {\it If $c = 12$, then $\pi_1(X^0)$ equals $\bold{Z}/(2)$
or $(\bold{Z}/(2))^{\oplus 2}$ (both are realizable; cf. Lemma $1.8$).}

\par
(2) {\it Suppose that $c \le 11$ and $\bold{Z}[\Delta]$ is non-primitive
in $H^2(\widetilde{X}, \bold{Z})$. Then $c \ge 8$ and
$\pi_1(X^0) = \bold{Z}/(2)$ (each $8 \le c \le 11$ is realizable).}

\par
(3) {\it If $c = 13$, then $\pi_1(X^0)$ equals $(\bold{Z}/(2))^{\oplus 2}$.}

\par
(4) {\it If $c = 14$, then $\pi_1(X^0)$ equals $(\bold{Z}/(2))^{\oplus 3}$.}

\par
(5) {\it If $c = 15$, then $\pi_1(X^0)$ equals $(\bold{Z}/(2))^{\oplus 4}$.}

\par \vskip 0.8pc
{\it Proof.} (1) By the proof of [Ni1, Lemma 4] or [B1, Lemma 3],
$\bold{Z}[\Delta]$ is not primitive. So there is a
double cover $\sigma : Z \rightarrow X$
ramified exactly at an 8-point subset $H$ of
$\text{\rm Sing} X$. One has
$\text{\rm Sing} Z = \sigma^{-1}(
\text{\rm Sing} X \setminus H)$,
consisting of 8 singular points of type $A_1$.
If $\text{\rm Sing} Z$ is not $2$-divisible,
then the condition in Lemma 1.9 is satisfied (Lemma 1.1),
whence $\pi_1(Z^0) = (1)$ and $\pi_1(X^0) = \bold{Z}/(2)$.
If $\text{\rm Sing} Z$ is $2$-divisible
then $\pi_1(X^0) = (\bold{Z}/(2))^{\oplus 2}$
as in Lemma 1.8.

\par
(2) follows from Lemma 1.3 and the arguments in (3).
For the realization of each $c$, we let $H$ be any affine hyperplane of
$\text{\rm Sing} A/\langle \iota \rangle$ (cf. 1.7)
and $X \rightarrow A/\langle \iota \rangle$ a minimal
resolution of any $16-c$ points not in $H$.

\par
(3) By Lemma 1.6, there are two 8-point subsets $H_1, H_2$
of $\text{\rm Sing} X$ with $|H_1 \cap H_2| = 4$.
Let $\sigma : Z \rightarrow X$
be the double cover ramified exactly over $H_1$.
Then $\widetilde{H}_2 := \sigma^{-1}(H_2 \setminus H_1)$
is a $2$-divisible 8-point subset of $\text{\rm Sing} Z$;
to see this, we apply 1.1, pull back the relation on $\widetilde{X}$
arising from the $2$-divisible set $H_2$ to
a relation on a minimal resolution of
$Z$ and apply 1.1 again. Note that $\text{\rm Sing} Z$
consists of 10 points of type $A_1$.
Let $\tau : Y \rightarrow Z$
be the double cover ramified exactly over $\widetilde{H}_2$.
Then $\text{\rm Sing} Y$
consists of 4 points of type $A_1$. So
$\pi_1(Y^0) = (1)$ by Lemma 1.9. Thus $|\pi_1(X^0)| = 4$
so that $Y/\pi_1(X^0) = X$.
Since $X$ has at worst type $A_1$ singularities,
$\pi_1(X^0) \cong (\bold{Z}/(2))^{\oplus 2}$
(cf. [X, Theorem 3]).

\par
(4) $c = 14$ implies that $\text{\rm Sing} X =
H_1 \cup H_2 \cup H_3$ with $|H_1 \cap H_2 \cap H_3| = 2$
(Lemmas 1.6 and 1.9). Let $\sigma : Z \rightarrow X$
be the double cover ramified exactly over $H_1$.
Set $\widetilde{H}_i = \sigma^{-1}(H_i \setminus H_1)$
($i = 2, 3$). Then $\text{\rm Sing} Z = \widetilde{H}_2
\cup \widetilde{H}_3$. As in (3), $\widetilde{H}_2$
and $\widetilde{H}_3$ are $2$-divisible. Hence $\pi_1(Z^0) =
(\bold{Z}/(2))^{\oplus 2}$ as in the proof of Lemma 1.8.
So $\pi_1(X^0) = (\bold{Z}/(2))^{\oplus 3}$ by the same
reasoning as in (3).

\par
(5) Let $H$ be a $2$-divisible 8-point subset of
$\text{\rm Sing} X$ (Lemma 1.9) and let
$\sigma : Z \rightarrow X$
be the double cover ramified exactly over $H$.
Then $\text{\rm Sing} Z$ consists of exactly
14 points of type $A_1$. Now (5) follows from (4)
and the reasoning in (3).

\par \vskip 2pc
{\bf \S 2. The  $K3$  case with $p \ge 3$}

\par \vskip 0.8pc
We shall prove Theorem 1 at the end of the section.
We treat first the case $p = 3$. Let us start with:

\par \vskip 0.8pc
{\bf Example 2.1.} For each $c \in \{1, \dots, 7\}$,
we shall construct an example of $X$
satisfying the conditions of Theorem 1 with $p = 3$
and $\pi_1(X^0) = (1)$; in particular, $\bold{Z}[\Delta]$
is primitive in $H^2(X, \bold{Z})$.

\par
It suffices to do for $c = 7$.
Let $\widetilde{X} \rightarrow \bold{P}^1$ be an elliptic $K3$  surface
with a section $P_0$, singular fibres of type
$I_1, I_1, I_2, I_3, I_4, I_{13}$
and trivial Mordell Weil group $MW$. This is No.39 in
[MP, the Table] or No.91 in [SZ, Table 2].
Clearly, $P_0$ together with some fibre components form
a divisor $\Delta$ of Dynkin type $7A_2$.
Let $\widetilde{X} \rightarrow X$ be the contraction of $\Delta$.
By [No, Lemma 1.5], if one lets $F$ be a general fibre,
then one has an exact sequence:
$$\pi_1(F \setminus P_0) \rightarrow \pi_1(X^0) \rightarrow \pi_1(\bold{P}^1)
= (1).$$
Note that the first homomorphism above factors through
$\pi_1(F_1 \setminus P_0)$ ($= \bold{Z}$) where $F_1$
is a fibre of type $I_1$. Hence $\pi_1(X^0)$ is cyclic.
Since the group $MW$ is trivial and all fibres are reducible,
the components of $\Delta$ form a partial $\bold{Z}$-basis
of $Pic(\widetilde{X})$ and hence $\bold{Z}[\Delta]$
is primitive in $H^2(\widetilde{X}, \bold{Z})$.
Thus the group $\pi_1(X^0)$ is perfect (Lemma 1.1);
so it is trivial.

\par \vskip 0.8pc
{\bf Example 2.2.} Here is an example of $X$
satisfying Theorem 1 with $(p, c) = (3, 8)$
and $\pi_1(X^0) = \bold{Z}/(3)$.

\par
Let $\widetilde{X} \rightarrow \bold{P}^1$ be an elliptic $K3$  surface
with a section $P_0$, singular fibres of type
$I_2, I_3, I_3, I_4, I_6, I_6$ and
the Mordell Weil group $MW \cong \bold{Z}/(6)$. This is No.108 in
[MP, the Table] or No.8 in [SZ, Table 2].
Write the 6 singular fibres as (in natural ordering)
$$\sum_{i=0}^1 G_i, \, \sum_{i=0}^2 A_i, \, \sum_{i=0}^2 B_i, \,
\sum_{i=0}^3 C_i, \, \sum_{i=0}^5 D_i, \, \sum_{i=0}^5 E_i,$$
so that $P_0$ meets components with index 0. Let $P_1$
be a generator of the group $MW$. By the height pairing in [Sh],
one can verify that (after relabeling) $P_1$ meets\
$G_0, A_1, B_1, C_2, D_1, E_1$ and $P_2 = 2P_1$ meets
$G_0, A_2, B_2, C_0, D_2, E_2$, and $P_1 \cap P_2 = \emptyset$.

\par
Let $\Delta = (P_0 + G_0) + (A_1+A_2) + (B_1+B_2) + (C_1+C_2) +
(D_1+D_2) + (D_4+D_5) + (E_1+E_2) + (E_4+E_5)$, which is of
Dynkin type $8A_2$. Expressing $P_2$ as a $\bold{Q}$-combination
of $P_0, F$ (a general fibre) and fibre components of index $\ge 1$,
we get:
$$\gather A_1+2A_2 + B_1+2B_2 + D_5+2D_4 +\\
D_2+2D_1 + E_5+2E_4 + E_2 + 2E_1
= 3L, \\
L = P_0 - P_2 + 2F - (D_3+D_2+E_3+E_4). \endgather$$
Let $f : \widetilde{X} \rightarrow X$ be the contraction of $\Delta$.
Denote by $\overline{L}$ the image on $X$ of $L$.
Then ${\Cal O}(\overline{L})^{\otimes 3} \cong \Cal{O}_X$.
Let $\sigma : Y = Spec \oplus_{i=0}^2 \Cal{O}_X(-i \, \overline{L})
\rightarrow X$
be the canonical Galois $\bold{Z}/(3)$-cover unramified
over $X^0$. Note that $Y$ consists of 6 points
of type $A_2$ (the preimages of $f(P_0+G_0)$, $f(C_1+C_2)$).

\par
Let $\widetilde{Y} \rightarrow Y$ be a minimal resolution
with $\Gamma$ the exceptional divisor.
Then the preimage on $\widetilde{Y}$ of $C_3$ is a disjoint
union of $C_3', C_3'', C_3'''$ so that $C_3' . \Gamma = 1$.
If $\text{\rm Sing} Y$ is $3$-divisible, then as in 1.1,
we get a relation:
$\sum_{i=1}^2 \sum_{k=1}^2 k \Gamma_i(k) \sim 3M$.
Intersecting this with $C_3'$ we see that $C_3'$
meets at least two components of $\Gamma$, a contradiction.
Hence $\text{\rm Sing} Y$ does not include any $3$-divisible subset
(cf. Lemma 1.3) and hence $\pi_1(Y^0) = (1)$ by Proposition 2.5 (1) below.
Thus $\pi_1(X^0) = \bold{Z}/(3)$.

\par \vskip 0.8pc
{\bf 2.3.} Let $A$ be an abelian surface with an order-3
symplectic automorphism
$\tau$ so that $A^{\tau}$ is a 9-point set. Such an example
is shown in [BL]. Then $X = A/\langle \tau \rangle$
is a $K3$  surface with 9 singularities of type $A_2$.

\par \vskip 0.8pc
{\bf Lemma 2.4.} {\it Let $X = A/\langle \tau \rangle$ be as in $2.3$.
In the following, we let $X_1 \rightarrow X$ be a minimal
resolution of singularities not in $\Sigma$
and $X_1^0 = X_1 \setminus \text{\rm Sing} X_1$.}

(1) {\it For each $c = 6, 7$, there is a $c$-point
subset $\Sigma$ of $\text{\rm Sing} X$ such that
$\pi_1(X_1^0) = \pi_1(X \setminus \Sigma) = \bold{Z}/(3)$.}

\par
(2) {\it There is an $8$-point
subset $\Sigma$ of $\text{\rm Sing} X$ such that
$\pi_1(X_1^0) = \pi_1(X \setminus \Sigma) = (\bold{Z}/(3))^{\oplus 2}$.}

\par \vskip 0.8pc
{\it Proof.} By [B1, Claim 2 in \S 4] : ``each pair of points
lie on a unique line'', which means that each 7-point
subset of $\text{\rm Sing} X$ includes a unique $3$-divisible
subset 6-point subset $H$. Let $\Sigma = H$ (resp. $\Sigma = H \cup \{p_1\}$
with $p_1$ a singular point of $X$ not in $H$) when
$c = 6$ (resp. $c = 7$).
Let $\sigma : Y \rightarrow X$ be the Galois $\bold{Z}/(3)$-cover
ramified exactly over $H$. Now $\hat{\Sigma} :=  \sigma^{-1}(\Sigma)
\cap \text{\rm Sing} Y$ consists of $3(c-6)$ points of
type $A_2$ and hence does not include any $3$-divisible subsets (Lemma 1.3).
So the group $\pi_1(Y \setminus \hat{\Sigma})$ is perfect (Lemma 1.1).
We have also $\pi_1(X \setminus \Sigma)/\pi_1(Y \setminus \hat{\Sigma})
\cong \bold{Z}/(3)$.
Now the group $\pi_1(Y \setminus \hat{\Sigma})$ is trivial because
it is also soluble being the subgroup of $\pi_1(X \setminus \Sigma)$,
while the latter is the image of the soluble group $\pi_1(X^0)$ [Remark 1.4].
This proves (1).

\par
(2) Let $\Sigma$ be an 8-point subset of $\text{\rm Sing} X$
including two 6-point subsets $H_1, H_2$ with $|H_1 \cap H_2| = 4$
(in notation of [B1], the lines determined by $H_1, H_2$
have a unique common point). Now (2) is similar to Lemma 1.8 (2)
or Proposition 1.10 (3).

\par \vskip 0.8pc
{\bf Proposition 2.5.} {\it Let $\widetilde{X} \rightarrow X, \Delta, c$
be as in Theorem $1$ with $p = 3$.}

\par
(1) {\it Suppose that $\bold{Z}[\Delta]$
is primitive in $H^2(X, \bold{Z})$, i.e., $\text{\rm Sing} X$
does not include any $3$-divisible subset (this is true if $c \le 5$).
Then $c \le 7$ and $\pi_1(X^0) = (1)$ (each $c \le 7$ is realizable
by Example $2.1$).}

\par
(2) {\it Suppose that $c \le 7$ and $\bold{Z}[\Delta]$
is non-primitive in $H^2(X, \bold{Z})$. Then $c = 6, 7$
and $\pi_1(X^0) = \bold{Z}/(3)$ (both $c$ are realizable
by Lemma $2.4$).}

\par
(3) {\it Suppose that $c = 8$. Then $\pi_1(X^0)$ equals
$\bold{Z}/(3)$, or $(\bold{Z}/(3))^{\oplus 2}$
(both groups are realizable by Example $2.2$ and Lemma $2.4$).}

\par \vskip 0.8pc
{\it Proof.} (1) As in [B1, Lemma 3], the
primitivity of $\bold{Z}[\Delta]$ implies that
$c \le 7$. Now as in Lemma 1.9, we are reduced
to show $\pi_1(X_1^0) = (1)$ for a particular $X_1$
satisfying Theorem 1 with $p = 3$ and $c \le 7$.
So just let $X_1$ be the one constructed in Example 2.1,
and (1) is proved.

\par
(2) By 1.1 and Lemma 1.3, one has $c \ge 6$, and
there is a $3$-divisible 6-point subset $H$ of
$\text{\rm Sing} X$ and a corresponding
Galois $\bold{Z}/(3)$-cover $Y \rightarrow X$
ramified exactly over $H$. Now $\text{\rm Sing} Y$
consists of $3(c-6)$ points of type $A_2$ and hence
does not include any $3$-divisible subsets (Lemma 1.3).
Thus $\pi_1(Y^0) = (1)$ by (1), whence $\pi_1(X^0) = \bold{Z}/(3)$.

\par
(3) This is similar to Proposition 1.10 (applying (1)).

\par \vskip 0.8pc
Next we consider the case $p > 3$. We begin with examples.

\par \vskip 0.8pc
{\bf Example 2.6.} (1) For each $c \le 4$, we construct an example
$X$ satisfying the conditions in Theorem 1
with $p = 5$ and $\pi_1(X^0) = (1)$; in particular,
${\bold Z}[\Delta]$ is primitive in $H^2(\widetilde{X}, {\bold Z})$.

\par
It suffices to construct an $X$ with $c = 4$.
Let $\widetilde{X} \rightarrow \bold{P}^1$ be an elliptic $K3$  surface
with a section $P_0$, singular fibres of type
$I_1, I_1, I_5, I_5, I_6, I_6$
and trivial Mordell Weil group $MW$. This is No.64 in
[MP, the Table] or No.9 in [SZ, Table 2].
Clearly, some fibre components form
a divisor $\Delta$ of Dynkin type $4A_4$.
Let $\widetilde{X} \rightarrow X$ be the contraction of $\Delta$.
Then as in Example 2.1, one has $\pi_1(X^0) = (1)$.

\par
(2) For each $c \le 3$, we construct an example
$X$ satisfying the conditions in Theorem 1
with $p = 7$ and $\pi_1(X^0) = (1)$; in particular,
${\bold Z}[\Delta]$ is primitive in $H^2(\widetilde{X}, {\bold Z})$.

\par
It suffices to construct an $X$ with $c = 3$.
Let $\widetilde{X} \rightarrow \bold{P}^1$ be an elliptic $K3$  surface
with a section $P_0$, singular fibres of type
$I_1, I_1, I_1, I_6, I_7, I_8$
and trivial Mordell Weil group. This is No.29 in
[MP, the Table] or No.41 in [SZ, Table 2].
Clearly, $P_0$ and some fibre components form
a divisor $\Delta$ of Dynkin type $3A_6$.
Let $\widetilde{X} \rightarrow X$ be the contraction of $\Delta$.
Then as in Example 2.1, one has $\pi_1(X^0) = (1)$.

\par \vskip 0.8pc
{\bf Example 2.7.} (1) We construct an example
$X$ satisfying the conditions in Theorem 1
with $(p, c) = (5, 4)$ and $\pi_1(X^0) = \bold{Z}/(5)$. Also see
Remark 3.3 for another construction.

\par
Let $\widetilde{X} \rightarrow \bold{P}^1$ be an elliptic $K3$  surface
with a section $P_0$, singular fibres of type
$I_1, I_1, I_1, I_1, I_{10}, I_{10}$ and
the Mordell Weil group $MW \cong \bold{Z}/(5)$. This is No.9 in
[MP, the Table] or No.54 in [SZ, Table 2].
Write the type $I_{10}$ singular fibres as (in natural ordering)
$$\sum_{i=0}^9 A_i, \,\, \sum_{i=0}^9 B_i$$
so that $P_0$ meets $A_0, B_0$. Then a generator $P_1$
of the group $MW$ meets $A_2, B_4$ after relabeling.
Let $\Delta = \sum_{i=1}^4 A_i + \sum_{i=6}^9 A_i +
\sum_{i=1}^4 B_i + \sum_{i=6}^9 B_i$.
Let $\widetilde{X} \rightarrow X$ be the contraction of $\Delta$.
As in Example 2.2,
we can verify that
$$\gather (4A_1+3A_2 +2A_3+A_4) + (4A_6+3A_7 +2A_8+A_9) +\\
(3B_1+B_2 +4B_3+2B_4) + (3B_6+B_7 +4B_8+2B_9)
= 5L, \\
L = P_0 - P_1 + 2F - (A_2+A_3+A_4+A_5+B_2+B_3+2B_4+2B_5+B_6+B_7),
\endgather$$
and proceed as there to obtain $\pi_1(X^0) = \bold{Z}/(5)$.

\par
(2) We construct an example $X$
satisfying the conditions in Theorem 1 with
$(p, c) = (7, 3)$ and $\pi_1(X^0) = \bold{Z}/(7)$.

\par
Let $\widetilde{X} \rightarrow \bold{P}^1$ be an elliptic $K3$  surface
with a section $P_0$, singular fibres of type
$I_1, I_1, I_1, I_7, I_7, I_7$ and
the Mordell Weil group $MW \cong \bold{Z}/(7)$. This is No.30 in
[MP, the Table] or No.13 in [SZ, Table 2].
Write the type $I_7$ singular fibres as (in natural ordering)
$$\sum_{i=0}^6 A_i, \,\, \sum_{i=0}^6 B_i, \,\, \sum_{i=0}^6 C_i$$
so that $P_0$ meets $A_0, B_0, C_0$. Then a generator $P_1$
of the group $MW$ meets $A_1, B_2, C_3$ after relabeling.
Let $\Delta = \sum_{i=1}^6 A_i + \sum_{i=1}^6 B_i +
\sum_{i=1}^6 C_i$.
Let $\widetilde{X} \rightarrow X$ be the contraction of $\Delta$.
As in Example 2.2, we can verify that
$$\gather (6A_1+5A_2+4A_3+3A_4+2A_5+A_6) + (5B_1+3B_2+B_3+6B_4+4B_5+2B_6) +\\
(4C_1+C_2+5C_3+2C_4+6C_5+3C_6) = 7L, \\
L = P_0 - P_1 + 2F - (B_2+B_3+C_2+C_3+C_4),
\endgather$$
and proceed as there to obtain $\pi_1(X^0) = \bold{Z}/(7)$.

\par \vskip 0.8pc
{\bf Example 2.8.} For each $k \le 18$, we construct a $K3$  surface $X$
so that $X$ has a type $A_k$ singularity as its only singularity and
$\pi_1(X^0) = (1)$.

\par
It suffices to construct an $X$ with $k = 18$.
Let $\widetilde{X} \rightarrow \bold{P}^1$ be an elliptic $K3$  surface
with a section $P_0$, singular fibres of type
$I_1, I_1, I_1, I_1, I_1, I_{19}$
and trivial Mordell Weil group. This is No.1 in
[MP, the Table] or No.112 in [SZ, Table 2].
Clearly, some fibre components form
a divisor $\Delta$ of Dynkin type $A_{18}$.
Let $\widetilde{X} \rightarrow X$ be the contraction of $\Delta$.
Then as in Example 2.1, one has $\pi_1(X^0) = (1)$.

\par \vskip 0.8pc
{\bf Proposition 2.9.} {\it Let $f : \widetilde{X} \rightarrow X$,
$\Delta$, $c \ge 1$, $p$ be as in Theorem $1$. Then we have:}

\par
(1) {\it $p \le 19$; if $p = 5$ then $c \le 4$; if $p = 7$ then
$c \le 3$; if $p > 7$ then $c = 1$.}

\par
(2) {\it Suppose that $p = 5, 7$ and $\bold{Z}[\Delta]$
is primitive in $H^2(\widetilde{X}, \bold{Z})$
(this is true when $p = 5$ and $c \le 3$, or $p = 7$ and $c \le 2$).
Then $\pi_1(X^0) = (1)$ (all $(p, c) = (5, \le 4), (7, \le 3)$
are realizable by Example $2.6$).}

\par
(3) {\it Suppose that $p = 5, 7$ and $\bold{Z}[\Delta]$
is non-primitive in $H^2(\widetilde{X}, \bold{Z})$.
Then $(p, c)$ equals $(5, 4)$ or $(7, 3)$, and $\pi_1(X^0) = \bold{Z}/(p)$
(both cases are realizable by Example $2.7$).}

\par
(4) {\it Suppose that $p > 7$. Then $\pi_1(X^0) = 1$
(all prime numbers $7 < p \le 19$ are realizable by Example $2.8$).}

\par \vskip 0.8pc
{\it Proof.} (1) follows from the calculation
$20 \ge \rho(\widetilde{X}) = \rho(X) + c(p-1) \ge 1 + c(p-1)$.

\par
As in Lemma 1.9, the assertions (2) and (4) need to be
verified only for a particular $X$ in Example 2.6 or 2.8,
and hence are true.

\par
(3) By Lemma 1.1, Sing $X$ is $p$-divisible. So
(3) follows from Lemma 1.3.

\par \vskip 0.8pc
Now Theorem 1 in the introduction is a consequence of Remark 1.4, Lemma 1.9 and
Propositions 1.10, 2.5 and 2.9.

\par \vskip 2pc
{\bf \S 3. The fundamental group of an open Enriques surface}

\par \vskip 0.8pc
We shall prove Theorem 2 in the section.
Let $W$ be an Enriques surface. The second cohomology
group $H^2(W,\bold{Z})\cong \text{\rm Pic} W$ is isomorphic to
$\bold{Z}^{10}\oplus \bold{Z}/(2)$, where the torsion is the
canonical class $K_W$. The free part $H^2(W,\bold{Z})_0$ admits
a canonical structure of a lattice which is even, unimodular and
of signature $(1,9)$ and hence
isomorphic to $U\oplus E_8$, where $U$ is the unimodular
hyperbolic lattice of signature $(1,1)$, and $E_8$ the negative
definite lattice associated with the Dynkin diagram of type $E_8$.

\par
Assume that $W$ contains a configuration of smooth rational curves of
Dynkin type $cA_{p-1}$, where $p$ is a prime. Then the pair $(p,c)$
is one of the following:

\par
$p=7, \quad c=1$

\par
$p=5, \quad c=1, 2$

\par
$p=3, \quad c=1, 2, 3, 4$

\par
$p=2, \quad c=1, 2, \cdots , 8.$

\par
Conversely, for each pair $(p,c)$ in the above list, by considering
various ellipic fibrations one can prove the existence of
an Enriques surface with $c$ singularities of type $A_{p-1}$
(see [CD]).

\par \vskip 0.8pc
Suppose that an Enriques surface  $W$ contains a configuration of
rational curves of Dynkin type $cA_{p-1}$. We fix the following
notation:

\par
$W^0$ = the open Enriques surface obtained by deleting those
$c(p-1)$ rational curves from $W$.

\par
$X^0$ = the inverse of $W^0$ in the $K3$  cover ${\widetilde X}$ of $W$.

\par \vskip 0.8pc
{\bf Lemma 3.1.} {\it If $(p,c)=(7,1)$, then
$\pi_1(W^0) =  \bold{Z}/(2)$. }

\par \vskip 0.8pc
{\it Proof.} In this case $X^0$ corresponds to the case
$(p,c)=(7,2)$ in Table 1, so it is simply connected.

\par \vskip 0.8pc
{\bf Lemma 3.2.} (1) {\it If $(p,c)=(5,1)$, then
$\pi_1(W^0) =  \bold{Z}/(2)$. }

\par
(2) {\it  If $(p,c)=(5,2)$, then $\pi_1(W^0) = \bold{Z}/(2)$,
$\bold{Z}/(10)$, or the dihedral group $D_{10}$ of order $10$. The
first case occurs if the $4A_4$ on the  $K3$  cover of $W$ is
primitive; the second if the $2A_4$ on $W$ is non-primitive; the
third if
the $2A_4$ on $W$ is primitive, while  the $4A_4$ on the $K3$  cover
is non-primitive. All three cases occur. (See  Examples $3.4$ below.)}

\par \vskip 0.8pc
{\it Proof.} (1) This case follows immediately from Table 1.

\par
(2) Since $\pi_1(W^0)$ is an extension of $\pi_1(X^0)$ by
$\bold{Z}/(2)$, we see from Table 1 that
$\pi_1(W^0)=\bold{Z}/(2)$, $\bold{Z}/(10)$, or the dihedral group
$D_{10}$ of order 10. The second group contains a normal subgroup
of index 5, and occurs as $\pi_1(W^0)$ only if there is a Galois
covering of $W$ of degree 5, unramified over $W^0$. The third
group contains no normal subgroup of index 5.

\par \vskip 0.8pc
{\bf Lemma 3.3.} {\it  Let $p$ be an odd prime. Let $D$ be a
divisor on an Enriques surface $W$. Suppose that {\rm Pic}$W$
contains a subgroup $N$ of finite index coprime to $p$, and that
the intersection number of $D$ with any element of $N$ is a
multiple of $p$. Then $D$ is $p$-divisible in {\rm Pic}$W$.}

\par \vskip 0.8pc
{\it Proof.} This follows from the unimodularity of
Pic$W$/(torsion) and the $p$-divisibility of the 2-torsion
$K_W=pK_W$.

\par \vskip 0.8pc
{\bf Examples 3.4.}

\par
{\bf (3.4.1)} The case with $(p,c)=(5,2)$ and $\pi_1(W^0) =
\bold{Z}/(10)$

\par
Let $W$ be the Example IV from [Kon]; this is one of the 7
families of Enriques surfaces with finite automorphisms.

\par
There are 20 smooth rational curves $E_1, ..., E_{20}$ on $W$.
(See Figure 4.4 in [Kon].) Take 8 curves $E_{16}$, $E_{4}$, $E_{3}$,
$E_{13}$, $E_{11}$, $E_{5}$, $E_{8}$, $E_{10}$ on $W$, which form
a configuration of Dynkin type $2A_4$. These are irreducible
components of an elliptic pencil of type $2I_5\oplus 2I_5$. We
claim that the divisor
$$D=(E_{16}+2E_{4}+3E_{3}+4E_{13})+(2E_{11}+4E_{5}+E_{8}+3E_{10})$$
is 5-divisible in Pic$W$. To see this, first note that $D$
intersects with any of the 20 curves $E_i$ in a multiple of 5
points. Next, consider an elliptic pencil of type $I_0^*\oplus
I_0^*$ together with a double section to infer that among the 20
curves are there 10 curves which generate a sublattice isomorphic
to
$$D_{4}\oplus D_4\oplus \pmatrix 0&2\\2&-2\endpmatrix,$$
a sublattice of index $2^3$ of the unimodular lattice
Pic$W$/(torsion). Now apply Lemma 3.3.

\par \vskip 0.8pc
{\bf (3.4.2)} The case with $(p,c)=(5,2)$ and $\pi_1(W^0) =
\bold{Z}/(2)$

\par
Let $W$ be the same surface as in Example (3.4.1).
 Take 8 curves $E_{16}$, $E_{4}$, $E_{3}$,
$E_{13}$, $E_{17}$, $E_{11}$, $E_{5}$, $E_{8}$ on $W$, which form
a configuration of Dynkin type $2A_4$. These are irreducible
components of the same elliptic pencil of type $2I_5\oplus 2I_5$
as above. The corresponding 16 curves on the  $K3$ -cover of $W$
form a configuration of Dynkin type $4A_4$, and can be found in
Figure 4.3 in [Kon]. It is checked that for any mod 5 nontrivial
integral linear combination of the 16 curves can one find a
smooth rational curve which intersects the combination in a
non-multiple of 5 points. So, the $4A_4$ is primitive.

\par \vskip 0.8pc
{\bf (3.4.3)} The case with $(p,c)=(5,2)$ and $\pi_1(W^0) =
D_{10}$

\par
Let $W$ be the Example I from [Kon] (see also [D]). There are 12
smooth rational curves $F_1, ..., F_{12}$ on $W$. We give the
dual graph below for the readers' convenience.

\par \vskip 2pc
\centerline{Figure 3.4.3} \vskip 2pc

Take 8 curves $F_9$, $F_3$, $F_4$, $F_5$, $F_1$, $F_8$, $F_7$,
$F_{10}$ on $W$, which form a configuration of Dynkin type
$2A_4$. These are irreducible components of an elliptic pencil of
type $II^*$. By intersecting with $F_{11}$ and $F_{12}$, we see
easily that the $2A_4$ is primitive.

\par
On the other hand, the corresponding 16 curves $F_i^{\pm}$ on the
 $K3$ -cover of $W$ form a 5-divisible configuration of Dynkin type
$4A_4$. To see this, note that the 16 curves are irreducible
components of an elliptic pencil of type $II^*\oplus II^*$, so
that
 the divisor
$$F_9^++2F_3^++3F_4^++4F_5^++2F_1^++4F_8^++F_7^++3F_{10}^+$$
 $$+4F_9^-+3F_3^-+2F_4^-+F_5^-+3F_1^-+F_8^-+4F_7^-+2F_{10}^-=5L,$$
$$L=F_9^-+F_3^-+F_4^-+F_5^-+F_6^-+F_1^-+F_8^-+2F_7^-+F_{10}^--F_6^+-F_7^+,$$
 is clearly
5-divisible.

\par \vskip 0.8pc
{\bf Lemma 3.5.} {\it Let $W$ be an Enriques surface.}

\par
(1) {\it If $W$ has an elliptic pencil with a singular fibre of
type $IV^*$, then the configuration of Dynkin type $3A_2$
consisting of non-central components of this fibre is primitive,
while the corresponding $12$ curves on the $K3$-cover of $W$ form
a $3$-divisible configuration of Dynkin type $6A_2$. }

\par
(2) {\it Any configuration of smooth rational curves on $W$ of
Dynkin type $4A_2$ contains exactly one $3$-divisible
sub-configuration of Dynkin type $3A_2$.}

\par \vskip 0.8pc
{\it Proof.} (1) Write the singular fibre as
$$F_1+2F_2+F_3+2F_4+F_5+2F_6+3F_7.$$
The 6 curves $F_1$, ..., $F_6$ form a non-primitive $3A_2$ if and
only if the divisor
$$F_1+2F_2+F_3+2F_4+F_5+2F_6$$
is $3$-divisible, if and only if a general fibre is $3$-divisible,
which is impossible, because no elliptic pencil on an Enriques
surface has a triple fibre.

\par
On the other hand, the corresponding 12 curves $F_{i}^{\pm}$ on
the  $K3$ -cover form a $3$-divisible configuration of Dynkin
type $6A_2$, as the 12 curves are irreducible components of an
elliptic pencil of type $IV^*\oplus IV^*$, and hence
 the divisor on the $K3$ -cover
$$F_1^++2F_2^++F_3^++2F_4^++F_5^++2F_6^+$$
 $$+2F_1^-+F_2^-+2F_3^-+F_4^-+2F_5^-+F_6^-=3L,$$
$$L=F_1^-+F_2^-+F_3^-+F_4^-+F_5^-+F_6^-+F_7^--F_7^+,$$
is clearly $3$-divisible.

\par
(2) Let $M$ be the sublattice of the unimodular lattice
Pic$W$/(torsion) generated by the given 8 curves of Dynkin type
$4A_2$. Let $\overline{M}$ be its primitive closure. Since the
discriminant group of $M$ is a 3-elementary group with 4
generators and the orthogonal complement $M^{\perp}$ has rank 2,
$\overline{M}/M$ must have order 3 or $3^2$. In other words, $M$
is not primitive and contains exactly one or four $3$-divisible
sub-configurations of Dynkin type $3A_2$. The second possibility
can be ruled out by the following claim and (1).

\par \vskip 0.8pc
{\bf Claim 3.5.1.} Any configuration of smooth rational curves on $W$
of Dynkin type $4A_2$ is equivalent, by a composition of
reflections in a smooth rational curve, to a configuration of the
same type consisting of irreducible components of an elliptic
pencil of type $IV^*\oplus I_3$ or $IV^*\oplus 2I_3$.

\par \vskip 0.8pc
To prove the claim, observe that $-$det$M^{\perp}$ is a perfect square,
so that we can find an isotropic element of $M^{\perp}$,
and hence a primitive isotropic element $A$ of Pic$W$ which is
orthogonal to the 8 curves. The divisor $A$ consists of an
elliptic configuration $B$ and, possibly, trees of smooth
rational curves, say, $E_i$. These trees may contain some of the
8 curves. Let $g$ be the composition of reflections in a smooth
rational curve $E_i$ which maps $A$ to $B$. Then $g$ maps the 8
curves to 8 smooth rational curves which are irreducible
components of the elliptic pencil $|2B|$. (A reflection is, in
general, not even an effective isometry, but in our case $g$ has
the desired property.) Finally, It is easy to check that if an
elliptic pencil on an Enriques surface contains 8 smooth rational
curves of Dynkin type $4A_2$, then it must be of type $IV^*\oplus
I_3$ or $IV^*\oplus 2I_3$.

\par \vskip 0.8pc
{\bf Lemma 3.6.} (1) {\it  If $(p,c)=(3,1)$, or $(3,2)$, then
$\pi_1(W^0) =  \bold{Z}/(2)$. }

\par
(2) {\it If $(p,c)=(3,3)$, then $\pi_1(W^0) =  \bold{Z}/(2)$,
$\bold{Z}/(6)$, or the symmetry group $S_3$ of order $6$. The
first case occurs if the $6A_2$ on the $K3$  cover of $W$ is
primitive; the second if the $3A_2$ on $W$ is non-primitive; the
third if
the $3A_2$ on $W$ is primitive, while  the $6A_2$ on the $K3$  cover
is non-primitive. All three cases occur. (See  Examples $3.7.1-3$ below.)}

\par
(3) {\it If $(p,c)=(3,4)$, then $\pi_1(W^0) = \bold{Z}/(6)$, or
$S_3\times \bold{Z}/(3)$. The first case occurs if the $8A_2$ on
the $K3$  cover of $W$ contains only one $3$-divisible $6A_2$;
the second if the $8A_2$ on the $K3$  cover is a union of two
$3$-divisible $6A_2$. The second case is supported by an example.
(See Example $3.7.4$ below.)}

\par \vskip 0.8pc
{\it Proof.} (1) These two cases follow from Table 1.

\par
(2) From Table 1, we see that
$\pi_1(W^0)$ is an extension of $(1)$ or
$\bold{Z}/(3)$ by $\bold{Z}/(2)$ and hence is isomorphic to
$\bold{Z}/(2)$, $\bold{Z}/(6)$, or $S_3$.

\par
(3) From Table 1, we see that
$\pi_1(W^0)$ is an extension of
$\pi_1(X^0)=\bold{Z}/(3)$, or $(\bold{Z}/(3))^{\oplus 2}$,
by $\bold{Z}/(2)$. There are 5 possibilities:
$\bold{Z}/(6)$, $S_3$, $(\bold{Z}/(3))^{\oplus 2} \times
\bold{Z}/(2)$, $S_3\times \bold{Z}/(3)$, or $G_{18/5}$, where the
last group is the nonabelian group of order $18$,
$$G_{18/5} = \langle a,b,c|a^3=b^3=c^2=1, ab=ba, aca=bcb=c \rangle.$$
By Lemma 3.5(2), the $4A_2$ is non-primitive, so $\pi_1(W^0)$ has
a normal subgroup of index 3. This rules out the second and fifth
possibilities.  Note that the third group has 4 normal subgroups of
index 3. The third case occurs if and only if the $4A_2$ on $W$
contains four different $3$-divisible $3A_2$, if and only if the
$4A_2$ on $W$ is of index $3^2$ in its primitive closure in
Pic$W$/(torsion). This is impossible again by Lemma 3.5(2).

\par \vskip 0.8pc
{\bf Example 3.7.}

\par
{\bf (3.7.1)} The case with $(p,c)=(3,3)$ and
$\pi_1(W^0)=\bold{Z}/(2)$.

\par
Let $W$ be the Example II from [Kon]. There are 12 smooth
rational curves $F_1, ..., F_{12}$ on $W$. We give the dual graph
below for the readers' convenience.

\par \vskip 2pc
\centerline{Figure 3.7.1} \vskip 2pc

Take 6 curves $F_1, F_2, F_5, F_6, F_9, F_{10}$ on $W$, which form a
configuration of Dynkin type $3A_2$, and let $W^0$ be the
surface with these 6 curves removed from $W$.
On the $K3$  cover of $W$ we have
12 curves
$$F_1^+, F_2^+, F_1^-, F_2^-, F_5^+, F_6^+, F_5^-, F_6^-, F_9^+, F_{10}^+,
F_9^-, F_{10}^-$$ which form a configuration of Dynkin type
$6A_2$. We claim that this $6A_2$ is primitive, whence
$\pi_1(W^0) = {\bold Z}/(2)$ by Lemma 3.6. Suppose that there is
an integral linear combination of the 12 curves
$$D=a_1^+F_1^+ +a_2^+F_2^+ +a_1^-F_1^- +a_2^-F_2^- + \cdots +
a_9^-F_9^- +a_{10}^- F_{10}^-$$
which is $3$-divisible in the Picard lattice of the $K3$  cover. Intersecting
$D$ with  $F_1^+$ and $F_4^+$,  we see that modulo 3
$$<D, F_1^+>=-2a_1^+ +a_2^+ \equiv 0, \quad <D, F_4^+>=a_1^+\equiv 0.$$
Thus $a_1^+ \equiv a_2^+ \equiv 0$. Similarly, intersecting $D$ with
$F_1^-, F_4^-, F_5^{\pm}, F_8^{\pm}, F_9^{\pm}, F_{12}^{\pm}$, we see that
all coefficients of $D$ are 0 modulo 3. This proves the claim.

\par \vskip 0.8pc
{\bf (3.7.2)} The case with $(p,c)=(3,3)$ and
$\pi_1(W^0)=\bold{Z}/(6)$.

\par
Let $W$ be the Example V from [Kon]. There are 20 smooth rational
curves $E_1, ..., E_{20}$ on $W$; see Figure 5.5 in [Kon].  Take 6
curves $E_{1}$, $E_{6}$, $E_{7}$, $E_{8}$, $E_{14}$, $E_{16}$, on
$W$, which form a configuration of Dynkin type $3A_2$. We claim
that the divisor
$$D=2E_{1}+E_{6}+2E_{7}+E_{8}+E_{14}+2E_{16}$$
is $3$-divisible in Pic$W$. To see this, first note that $D$
intersects with any of the 20 curves $E_i$ in a multiple of 3
points. Next, consider the elliptic pencil $|E_{16}+E_{20}|$,
which is of type $III*\oplus 2I_2$. Its irreducible components
together with a double section $E_{18}$ generate a sublattice
isomorphic to
$$E_{7}\oplus A_1\oplus \pmatrix 0&1\\1&-2\endpmatrix,$$
a sublattice of index $2$ of the unimodular lattice
Pic$W$/(torsion). Now apply Lemma 3.3.

\par \vskip 0.8pc
{\bf (3.7.3)} The case with $(p,c)=(3,3)$ and $\pi_1(W^0)=S_3$.

\par
Let $W$ be an Enriques surface with an elliptic pencil containing
a singular fibre of type $IV^*$. Take the 6 curves of Dynkin type
$3A_2$ out of this fibre. Then the result follows from Lemma
3.5(1).

\par \vskip 0.8pc
{\bf (3.7.4)} The case with $(p,c)=(3,4)$ and
$\pi_1(W^0)=S_3\times \bold{Z}/(3)$.

\par
Let $W$ be the Example V from [Kon]. There are 20 smooth rational
curves $E_1, ..., E_{20}$ on $W$; see Figure 5.5 in [Kon]. Take 8
curves $E_{3}$, $E_{4}$, $E_{1}$, $E_{6}$, $E_{7}$, $E_{8}$,
$E_{14}$, $E_{16}$, on $W$, which form a configuration of Dynkin
type $4A_2$. These are irreducible components of an elliptic
pencil of type $IV^*\oplus I_3$. We have proved in Example 3.7.2
that the divisor
$$D=2E_{1}+E_{6}+2E_{7}+E_{8}+E_{14}+2E_{16}$$
is $3$-divisible in Pic$W$. On the other hand, by Lemma 3.5(1),
the 6 curves $E_{3}$, $E_{4}$, $E_{1}$, $E_{6}$, $E_{7}$, $E_{8}$
form a primitive configuration of Dynkin type $3A_2$, whose pull
back on the $K3$ -cover form a $3$-divisible configuration of
Dynkin type $6A_2$.

\par \vskip 0.8pc
{\bf Definition 3.8.} Let $W$ be an Enriques surface with a configuration
of Dynkin type $kA_1$, i.e. mutually disjoint $k$ smooth rational curves.
The configuration is called {\it $2$-divisible $k$-point set} if the sum
of the $k$ curves is equal to $2L$ for an integral divisor
$L$ on $W$; since $K_W$ is the only torsion element
in Pic($W$) and since $2L = 2(L + K_W)$, there are exactly
two double covers of $W$ both branched exactly at these
$k$ curves.

\par
Let ${\widetilde X}$ be the  $K3$  cover of $W$.
Then the pull back on ${\widetilde X}$ of a Dynkin type
$cA_1$ configuration on $W$, is of Dynkin type $2cA_1$.
Hence a configuration of Dynkin type $kA_1$ is $2$-divisible
only if $k=4$, or 8.
Note also that the pull back on ${\widetilde X}$
of $4A_1$ on $W$ is
$2$-divisible if and only if the $4A_1$ is congruent to 0
or $K_W$ modulo 2 in Pic($W$).

\par
Let $K_1$ and $K_2$ be distinct $2$-divisible $4$-point sets
on an Enriques surface.
Then $|K_1\cap K_2|=0$, or 2. If  $|K_1\cap K_2|=2$, then the symmetric
difference  $K_1\triangle K_2$ is also a $2$-divisible 4-point set.

\par \vskip 0.8pc
{\bf Lemma 3.9.} (1) {\it If $p=2$, $c=1,2$, or $3$, then
$\pi_1(W^0) =  \bold{Z}/(2)$. }

\par
(2) {\it If $(p,c)=(2,4)$, then
$\pi_1(W^0) =  \bold{Z}/(2)$, $(\bold{Z}/(2))^{\oplus 2}$, or $\bold{Z}/(4)$.
The first case occurs if the $8A_1$ on the  $K3$  cover of $W$ is
primitive; the second if the $4A_1$ on $W$ is $2$-divisible; the third if
 the $4A_1$ on $W$ is primitive, while  the $8A_1$ on the  $K3$  cover
 is $2$-divisible.}

\par
(3) {\it If $(p,c)=(2,5)$, then
$\pi_1(W^0) =  \bold{Z}/(2)$, $(\bold{Z}/(2))^{\oplus 2}$, or $\bold{Z}/(4)$.
The first case occurs if the $10A_1$ on the  $K3$  cover of $W$ is
primitive; the second if the $5A_1$ on $W$ contains a $2$-divisible $4$-point
subset; the third if
 the $5A_1$ on $W$ is primitive, while  the $10A_1$ on the  $K3$  cover
contains a $2$-divisible $8$-point subset.}

\par
(4) {\it If $(p,c)=(2,6)$, then
$\pi_1(W^0) =  \bold{Z}/(4)$, $(\bold{Z}/(2))^{\oplus 2}$,
$\bold{Z}/(4)\times \bold{Z}/(2)$, or $(\bold{Z}/(2))^{\oplus 3}$.
The first case occurs if the $6A_1$ on $W$ contains no
$2$-divisible $4$-point subset; the second if
the $12A_1$ on the  $K3$  cover of $W$ contains only
one $2$-divisible $8$-point subset
and the $6A_1$ on $W$ contains a
$2$-divisible $4$-point subset; the third if the $12A_1$ on the  $K3$  cover is
a union of two $2$-divisible $8$-point
subsets and the $6A_1$ on $W$ contains
only one $2$-divisible $4$-point subset; the fourth if the $6A_1$ on $W$ is
a union of two $2$-divisible $4$-point
subsets.}

\par
(5) {\it If $(p,c)=(2,7)$, then
$\pi_1(W^0) =(\bold{Z}/(2))^{\oplus 4}$,
$(\bold{Z}/(2))^{\oplus 2}\times \bold{Z}/(4)$, or $\Gamma_2c_1$, where
$$\Gamma_2c_1=<a,b,c|a^4=b^2=c^2=1, ab=ba, ac=ca^3b, bc=cb>.$$
The first case occurs if the $7A_1$ on $W$ is a union of three
$2$-divisible $4$-point subsets; the second if the $7A_1$ on $W$ is a union of
one $A_1$ and two
$2$-divisible $4$-point subsets; the third if the $7A_1$ on $W$ contains only
one $2$-divisible $4$-point subset. }

\par
(6) {\it If $(p,c)=(2,8)$, then
$\pi_1(W^0) =  (\bold{Z}^{\oplus 2}\rtimes \bold{Z}/(2))\rtimes \bold{Z}/(2)$. }

\par \vskip 0.8pc
All cases are supported by examples except the case with
 $(p,c)=(2,7)$ and
$\pi_1(W^0) =(\bold{Z}/(2))^{\oplus 4}$. (See Examples 3.12.)

\par \vskip 0.8pc
{\it Proof.} (1) and (6) follow from Table 1.

\par
(2) and (3) also follow from Table 1. Note that if a subconfiguration of Dynkin
type $4A_1$ on $W$ is $2$-divisible, i.e. $4A_1$ is linearly equivalent to
$2L$ for some $L\in$ Pic($W$), then both $L$ and $L+K_W$ determine Galois
double covers of $W$, which correspond to two of the three normal subgroups
of $(\bold{Z}/(2))^{\oplus 2}$ of index 2.

\par
(4) From Table 1, we see that
$\pi_1(W^0)$ is an extension of
$\bold{Z}/(2)$ by ${\bold Z}/(2)$ or
$(\bold{Z}/(2))^{\oplus 2}$. There are 5
possibilities: $\bold{Z}/(4)$, $(\bold{Z}/(2))^{\oplus 2}$,
$\bold{Z}/(4)\times \bold{Z}/(2)$, $(\bold{Z}/(2))^{\oplus 3}$, or
the dihedral group $D_8$ of order 8. The last can be ruled out by observing that
if $\pi_1(X^0)=(\bold{Z}/(2))^{\oplus 2}$, then $\pi_1(W^0)$ must have
an odd number of normal subgroups isomorphic to $(\bold{Z}/(2))^{\oplus 2}$, while
$D_8$ has exactly two such subgroups.

\par
(5) In this case,
$\pi_1(W^0)$ is an extension of ${\bold Z}/(2)$ by
$(\bold{Z}/(2))^{\oplus 3}$. There are 4
possibilities: $(\bold{Z}/(2))^{\oplus 4}$,
$\bold{Z}/(4)\times (\bold{Z}/(2))^{\oplus 2}$, $\Gamma_2c_1$, or
$D_8\times \bold{Z}/(2)$. The last can be ruled out by noting that
a double cover of $W$ branched along the union of four out of the seven
 curves is again an Enriques surface (with 4 points blown up) and that, by (4),
no open Enriques surface with $(p,c)=(2,6)$ can have $D_8$ as its fundamental
group.

\par
Note that
$\bold{Z}/(4)\times (\bold{Z}/(2))^{\oplus 2}$(resp. $\Gamma_2c_1$)
has exactly 7 (resp. 3) normal subgroups of index 2.

\par
The group $(\bold{Z}/(2))^{\oplus 4}$ has 15 normal subgroups of index 2,
and hence occurs as $\pi_1(W^0)$ only if the $7A_1$ contains 7 different
 $2$-divisible 4-point subsets. This condition is equivalent to that
the $7A_1$ is a union $K_1\cup K_2\cup K_3$ of three  $2$-divisible 4-point
subsets $K_i$, $i=1,2,3$, where the seven $2$-divisible 4-point subsets are
$K_1, K_2, K_3, K_1\triangle K_2, K_2\triangle K_3, K_1\triangle K_3$ and
$K_1\triangle K_2\triangle K_3$.

\par \vskip 0.8pc
{\bf 3.10.} Let $\widetilde{X}$ be the Kummer surface $Km(E_1\times E_2)$, where
$E_i$ is an elliptic curve with fundamental period $\tau_i$. Let $(a_1,a_2)$
be the 2-torsion point $((1+\tau_1)/2,(1+\tau_2)/2)\in E_1\times E_2$ and
consider the following involution of $E_1\times E_2$
$$\sigma : (z_1,z_2) \to (-z_1+a_1, z_2+a_2).$$
Then $\sigma$ induces a fixed point free involution $\bar{\sigma}$ on
 $\widetilde{X}$ and the quotient surface
$W_{\tau_1,\tau_2}= \widetilde{X}/\bar{\sigma}$ is
an Enriques surface. On $W_{\tau_1,\tau_2}$ we have 12 smooth rational
curves coming from
the sixteen 2-torsion points, (a 2-torsion) $\times E_2$, and
$E_1\times$ (a 2-torsion). Their dual graph is given in Figure 3.10.

\par \vskip 2pc
\centerline{Figure 3.10}
\par \vskip 2pc

There contained in the graph are
16 different configurations of type $I_8$,
half of them giving elliptic pencils on $W_{\tau_1,\tau_2}$
and the other half corresponding
to half elliptic pencils. We may assume that
$|F_1+F_2+F_3+F_4+F_5+F_6+F_7+F_8|$ is an elliptic pencil. Then
$|2(F_1+F_2+F_3+F_4+F_5+F_6+F_7+F_{10})|$ is also an elliptic pencil.
Modulo 2 in Pic($W_{\tau_1,\tau_2}$) there are many congruences.
To raise a few, we have
the following :

\par
(1) $F_1+F_2+F_3+F_4+F_5+F_6+F_7+F_8\equiv 0$ mod 2 in
 Pic($W_{\tau_1,\tau_2}$).

\par
(2) $F_2+F_4+F_9+F_{11}\equiv F_2+F_6+F_9+F_{12}\equiv 0$ mod 2.

\par
(3) $F_1+F_3+F_5+F_7\equiv F_2+F_4+F_6+F_8\equiv K_W$ mod 2.

\par \vskip 0.8pc
{\bf 3.11.}
If  $\tau_1=\tau_2=\sqrt{-1}$, then the special Enriques surface
$W_{\sqrt{-1},\sqrt{-1}}$ has additional 8 smooth rational curves,
$F_{13}, F_{14}, \cdots, F_{20}$ [Kon, Example III].
Their dual graph is [Kon, Fig 3.5, p. 212],
but we will use $F_i$ instead of $E_i$ in [Kon].

\par \vskip 0.8pc
{\bf Examples 3.12.}

\par
{\bf (3.12.1)} $(p,c)=(2,4)$ and $\pi_1(W^0) =  \bold{Z}/(2)$

\par
$W=W_{\tau_1,\tau_2}$.
Take the 4 curves, $F_2, F_4, F_6$, and $F_9$. Then the sum of the 8 curves
on  $\widetilde{X}$ has intersection number 1 with $F_7^+$.
Here we denote by $F_i^{+} \cup F_i^{-}$ the inverse
on ${\widetilde X}$ of $F_i$.

\par \vskip 0.8pc
{\bf (3.12.2)} $(p,c)=(2,4)$ and $\pi_1(W^0) =  (\bold{Z}/(2))^{\oplus 2}$

\par
$W=W_{\tau_1,\tau_2}$.
Take the 4 curves, $F_2, F_4, F_9$, and $F_{11}$. These form a $2$-divisible
4-point set (3.10.(2)).

\par \vskip 0.8pc
{\bf (3.12.3)} $(p,c)=(2,4)$ and $\pi_1(W^0) =  \bold{Z}/(4)$

\par
$W=W_{\tau_1,\tau_2}$.
Take the 4 curves, $F_2, F_4, F_6$, and $F_8$. Use 3.10.(3).

\par \vskip 0.8pc
{\bf (3.12.4)} $(p,c)=(2,5)$ and $\pi_1(W^0) =  \bold{Z}/(2)$

\par
$W=W_{\sqrt{-1},\sqrt{-1}}$.
Take the 5 curves, $F_4, F_6, F_8, F_9$, and $F_{14}$.

\par \vskip 0.8pc
{\bf (3.12.5)} $(p,c)=(2,5)$ and $\pi_1(W^0) =  (\bold{Z}/(2))^{\oplus 2}$

\par
$W=W_{\tau_1,\tau_2}$.
Take the 5 curves, $F_2, F_4, F_6, F_9$, and $F_{11}$.

\par \vskip 0.8pc
{\bf (3.12.6)} $(p,c)=(2,5)$ and $\pi_1(W^0) =  \bold{Z}/(4)$

\par
$W=W_{\tau_1,\tau_2}$.
Take the 5 curves, $F_2, F_4, F_6, F_8$, and $F_9$.

\par \vskip 0.8pc
{\bf (3.12.7)} $(p,c)=(2,6)$ and $\pi_1(W^0) =  \bold{Z}/(4)$

\par
$W=W_{\sqrt{-1},\sqrt{-1}}$.
Take the 6 curves, $F_2, F_4, F_6, F_8, F_{10}$, and $F_{16}$.

\par \vskip 0.8pc
{\bf (3.12.8)} $(p,c)=(2,6)$ and $\pi_1(W^0) =  (\bold{Z}/(2))^{\oplus 2}$

\par
$W=W_{\sqrt{-1},\sqrt{-1}}$.
Take the 6 curves, $F_4, F_6, F_8, F_{10}, F_{11}$, and $F_{15} $.

\par \vskip 0.8pc
{\bf (3.12.9)} $(p,c)=(2,6)$ and $\pi_1(W^0) =  \bold{Z}/(4)\times \bold{Z}/(2)$

\par
$W=W_{\tau_1,\tau_2}$.
Take the 6 curves, $F_4, F_6, F_8, F_9, F_{10}$, and $F_{12}$.

\par \vskip 0.8pc
{\bf (3.12.10)} $(p,c)=(2,6)$ and $\pi_1(W^0) =  (\bold{Z}/(2))^{\oplus 3}$

\par
$W=W_{\tau_1,\tau_2}$.
Take the 6 curves, $F_4, F_6, F_8, F_{10}, F_{11}$, and $F_{12}$.

\par \vskip 0.8pc
{\bf (3.12.11)} $(p,c)=(2,7)$ and $\pi_1(W^0) =
(\bold{Z}/(2))^{\oplus 2}\times \bold{Z}/(4)$

\par
$W=W_{\tau_1,\tau_2}$.
Take the 7 curves, $F_4, F_6, F_8, F_9, F_{10}, F_{11}$, and $F_{12}$.

\par \vskip 0.8pc
{\bf (3.12.12)} $(p,c)=(2,7)$ and $\pi_1(W^0) =\Gamma_2c_1$

\par
$W=W_{\sqrt{-1},\sqrt{-1}}$.
Take the 7 curves, $F_4, F_8, F_9, F_{10}, F_{11}, F_{12}$, and $F_{20}$.

\par \vskip 0.8pc
Combining results in this section, we conclude Theorem 2.

\par \vskip 2pc
{\bf \S 4. The proof of Theorem 3}

\par \vskip 0.8pc
We now prove Theorem 3.

\par \vskip 0.8pc

{\bf Claim 1.} H is either a $K3$  or an Enriques surface
with at worst Du Val singularities of type $A_{p-1}$.

\par \vskip 0.8pc
Note that $2K_H = 2(K_V + H)|H \sim 0$. So we have only to
show that $H^1(H, {\Cal O}_H) = 0$ (and hence the irregularity
of the resolution of $H$ also vanishes because $H$
has only rational singularities).
Consider the exact sequence:
$$0 \rightarrow {\Cal O}_V(-H) \rightarrow {\Cal O}_V
\rightarrow {\Cal O}_H \rightarrow 0.$$
This induces a long exact sequence of cohomologies.
Now the Kawamata-Viehweg vanishing theorem
implies that $H^1(V, {\Cal O}_V) = 0 = H^2(V, {\Cal O}(-H))$,
whence $H^1(X, {\Cal O}_H)$ $= 0$. This proves Claim 1.

\par \vskip 0.8pc
Embed $V$ in a projective space and let $L$ be a general
hyperplane on $V$ such that $L \cap$ Sing $H = \emptyset$
and $L \cap H$ is a smooth irreducible curve on $L$,
whence $L$ is smooth along this curve because $H$ is Cartier.
This is possible because the normal surface $H$ has only
finitely many singular points.
By the result of Hamm-Le in [HL, Theorem 1.1.3],
one has $\pi_1(V^0) = \pi_1(L \setminus L \cap$ Sing $V)$.

\par \vskip 0.8pc
{\bf Claim 2.} The natural homomorphism
$\pi_1(H \cap L \setminus H \cap L \cap $ Sing V)
$\rightarrow \pi_1(L \setminus L \cap $ Sing $V)$
is surjective.

\par \vskip 0.8pc
Let ${\widetilde L} \rightarrow L$ be the minimal resolution.
By the assumption, $H \cap L$ is away from Sing $L$, and hence
the pull back on ${\widetilde L}$,
denoted also by $H \cap L$, of $H \cap L$
is still smooth irreducible and also nef and big.
Note that $H \cap$ Sing $V \subseteq$ Sing $H$ because
$H$ is Cartier. Hence $H \cap L \cap$Sing $V = \emptyset$
by the choice of $L$; similarly, $\Delta :=
L \cap$ Sing $V \subseteq$ Sing $L$.
By [No, Cor. 2.3 and the proof of Cor. 2.4B],
we obtain the surjectivity of the homomorphism
$\pi_1(H \cap L) \rightarrow
\pi_1({\widetilde L} \setminus {\widetilde \Delta})
= \pi_1(L \setminus \Delta)$,
where ${\widetilde \Delta}$ is the inverse of $\Delta$
and the latter equality comes from the observation that
${\widetilde L} \setminus {\widetilde \Delta} \rightarrow L \setminus \Delta$
is the minimal resolution of singular points in (Sing $L) \setminus \Delta$
and the fact that every singular point on $L$ is log terminal because
so is $V$ and the generality of $L$ [Ko1, Theorem 7.8]. This proves Claim 2.

\par \vskip 0.8pc
Combining Claim 2 with the equality preceding it, we get a surjective homomorphism
$\pi_1(H \cap L \setminus H \cap L \cap $ Sing $V) \rightarrow
\pi_1(V^0)$. Since the above map factors
through $\pi_1(H \setminus H \cap$ Sing $V) \rightarrow
\pi_1(V^0)$, the latter map
is also surjective.

\par
On the other hand,
$H \cap$ Sing $V \subseteq$ Sing $H$, whence we have an
inclusion $H^0 := H \setminus$ Sing $H \subseteq H \setminus H \cap $Sing $V$
and its induced surjective homomorphism
$\pi_1(H^0) \rightarrow \pi_1(H \setminus H \cap $Sing $V)$.
This, combined with the early sujective map in the
preceding paragraph, produces a surjective homomorphism
$\pi_1(H^0) \rightarrow \pi_1(V^0)$. This, together
with Claim 1 and Theorems 1 and 2, implies Theorem 3.

\par \vskip 1pc
{\bf Added in proof.} After the paper was submitted, we learnt that
Conjecture B has been proved by S. Takayama under even weaker condition
[Ta], though Conjecture A is still open.

\par \vskip 2.5pc
{\bf References}

\par \vskip 0.3pc \noindent
[A] F. Ambro, Ladders on Fano varieties, alg-geom/{\bf 9710005}.

\par \vskip 0.3pc\noindent
[B1] W. Barth,  $K3$  surfaces with nine cusps, alg-geom/{\bf 9709031}.

\par \vskip 0.3pc \noindent
[B2] W. Barth, On the classification of  $K3$  surfaces with nine cusps,
math.AG/ {\bf 9805082}.

\par \vskip 0.3pc \noindent
[BL] C. Birkenhake and H. Lange, A family of abelian surfaces and curves
of genus four, Manuscr. Math. {\bf 85} (1994), 393--407.

\par \vskip 0.3pc \noindent
[C] F. Campana, On twistor spaces of the class $C$, J. Diff. Geom.
{\bf 33} (1991), 541--549.

\par \vskip 0.3pc \noindent
[CD] F. Cossec and I. Dolgachev, Enriques surfaces I,
Progress in Math. {\bf Vol. 76} (1989), Birkhauser.

\par \vskip 0.3pc \noindent
[D] I. Dolgachev, On automorphisms of Enriques surfaces, Invent. Math.
{\bf 76} (1984), 163-177.

\par \vskip 0.3pc \noindent
[FKL] A. Fujiki, R. Kobayashi and S. Lu, On the fundamental group of
certain open normal surfaces, Saitama Math. J. {\bf 11} (1993), 15--20.

\par \vskip 0.3pc \noindent
[GZ 1, 2] R. V. Gurjar and D. -Q. Zhang, $\pi_1$ of smooth points of
a log del Pezzo surface is finite, I; II, J. Math. Sci. Univ. Tokyo
{\bf 1} (1994), 137--180; {\bf 2} (1995), 165--196.

\par \vskip 0.3pc \noindent
[GM] M. Goresky and R. MacPherson, Stratified Morse theory, Springer,
1988.

\par \vskip 0.3pc \noindent
[HL] H. A. Hamm and D. T. Le, Lefschetz theorems on
quasi-projective varieties, Bull. Soc. math. France,
{\bf 113} (1985), 123-142.

\par \vskip 0.3pc \noindent
[Ka] Y. Kawamata, The cone of curves of algebraic varieties.
Ann. Math. {\bf 119} (1984), 603--633.

\par \vskip 0.3pc \noindent
[KM] S. Keel and J. McKernan, Rational curves on
quasi-projective varieties, Mem. Amer. Math.
Soc. {\bf 140} (1999).

\par \vskip 0.3pc \noindent
[Kol1] J. Kollar, Shafarevich maps and plurigenera of algebraic
varieties, Invent. Math. {\bf 113} (1993), 177--215.

\par \vskip 0.3pc \noindent
[Kol2] J. Kollar, Shafarevich maps and automorphic forms,
M. B. Porter Lectures at Rice Univ., Princeton Univ. Press.

\par \vskip 0.3pc \noindent
[Ke] J. Keum, Note on elliptic  $K3$  surfaces, Trans. A.M.S.
{\bf 352} (2000), 2077-2086.

\par \vskip 0.3pc \noindent
[KoMiMo] J. Kollar, Y. Miyaoka and S. Mori, Rationally connected varieties,
J. Alg. Geom. {\bf 1} (1992), 429--448.

\par \vskip 0.3pc \noindent
[Kon] S. Kondo, Enriques surfaces with finite automorphism groups,
Japan J. Math. {\bf 12} (1986), 191-282.

\par \vskip 0.3pc \noindent
[Mi] T. Minagawa, Deformations of {\bf Q}-Calabi-Yau 3-folds
and {\bf Q}-Fano 3-folds of Fano index 1, math.AG/{\bf 9905106}.

\par \vskip 0.3pc \noindent
[MP] R. Miranda and U. Persson, Mordell-Weil groups of extremal
elliptic  $K3$  surfaces, in :
Problems in the theory of surfaces and their classification
(Cortona, 1988), Sympo. Math. XXXII, Acad. Press, London,
1991, pp. 167--192.

\par \vskip 0.3pc \noindent
[Mu] S. Mukai, Finite groups of automorphisms of  $K3$  surfaces
and the Mathieu group, Invent. Math. {\bf 94} (1988), 183--221.

\par \vskip 0.3pc \noindent
[Ni1] V. V. Nikulin, On Kummer surfaces, Math. USSR Izv.
{\bf 9} (1975), 261--275.

\par \vskip 0.3pc \noindent
[Ni2] V. V. Nikulin, Integral symmetric bilinear forms
and some of their applications, Math. USSR Izv. {\bf 14} (1980),
103--167.

\par \vskip 0.3pc \noindent
[Ni3] V. V. Nikulin, Finite automorphism groups of Kahler  $K3$
surfaces, Trans. Moscow Math. Soc. {\bf 38} (1980), 71--135.

\par \vskip 0.3pc \noindent
[No] M. V. Nori, Zariski's conjecture and related problems,
Ann. Sci. Ecole Norm. Sup. {\bf 16} (1983), 305--344.

\par \vskip 0.3pc \noindent
[SZ] I. Shimada and D. -Q. Zhang, Classification of extremal elliptic
 $K3$  surfaces and fundamental groups of open  $K3$  surfaces,
Nagoya Math. J. to appear.

\par \vskip 0.3pc \noindent
[Sa] T. Sano, On classification of non-Gorenstein {\bf Q}-Fano
3-folds of Fano index 1, J. Math. Soc. Japan, {\bf 47} (1995), 369-380.

\par \vskip 0.3pc \noindent
[Sh] T. Shioda, On the Mordell-Weil lattices, Comment.
Math. Univ. Sancti pauli, {\bf 39} (1990), 211--240.

\par \vskip 0.3pc \noindent
[Ta] S. Takayama, Simple connectedness of weak Fano varieties,
J. Alg. Geom. {\bf 9} (2000), 403--407.

\par \vskip 0.3pc \noindent
[T] D. Toledo, Projective varieties with non-residually
finite fundamental group, Inst. Hautes
Etudes Sci. Publ. Math. {\bf 77} (1993), 103--119.

\par \vskip 0.3pc \noindent
[X] G. Xiao, Galois covers between  $K3$ surfaces,
Ann. Inst. Fourier (Grenoble), {\bf 46} (1996), 73--88.

\par \vskip 0.3pc \noindent
[Z1] D. -Q. Zhang, The fundamental group of the smooth
part of a log Fano variety, Osaka J. Math. {\bf 32} (1995), 637--644.

\par \vskip 0.3pc \noindent
[Z2] D. -Q. Zhang, Algebraic surfaces with nef and big anti-canonical divisor,
Math. Proc. Cambridge Philos. Soc. {\bf 117} (1995), 161--163.

\par \vskip 0.3pc \noindent
[Z3] D. -Q. Zhang, Algebraic surfaces with log canonical singularities and
the fundamental groups of their smooth parts, Trans. Amer. Math. Soc.
{\bf 348} (1996), 4175--4184.

\par \vskip 2.5pc \noindent
J. Keum
\par \noindent
Korea Institute for Advanced Study
\par \noindent
207-43 Cheongryangri-dong, Dongdaemun-gu
\par \noindent
Seoul 130-012, Korea
\par \noindent
e-mail : jhkeum$\@$kias.re.kr

\par \vskip 0.8pc \noindent
D. -Q. Zhang
\par \noindent
Department of Mathematics, National University of Singapore
\par \noindent
2 Science Drive 2, Singapore 117543
\par \noindent
Republic of Singapore
\par \noindent
e-mail : matzdq$\@$math.nus.edu.sg

\par \vskip 2pc \noindent
Figures 3.4.3, 3.7.1 and 3.10
are respectively Figures 1.4, 2.4 and 3.5 
(with $\sum_{i=1}^{12} E_i$ there replaced by
$\sum_{i=1}^{12} F_i$ here) in [Kon].

\par \vskip 2pc \noindent
In Table 1 below, we write $\text{\rm Sing} X = cA_{p-1}$, $X^0 = X - \text{\rm Sing} X$.
Let $f : {\widetilde X} \rightarrow X$ be the minimal resolution
with $D = f^{-1}(\text{\rm Sing} X)$. $H$ or $H_i$ is a 2-divisible configuration
of ${\bold P}^1$'s of Dynkin type $8A_1$.
$R$ or $R_i$ is a 3-divisible configuration of ${\bold P}^1$'s of Dynkin 
type $6A_2$.

\par \vskip 1pc \noindent
{\bf Table 1}

\par \vskip 1pc \noindent
No.1: $p = 2$; $1 \le c \le 11$; ${\bold Z}[D]$ is primitive in 
$H^2({\widetilde X}, {\bold Z})$; $\pi_1(X^0) = (1)$; $Y = X$.

\par \noindent
No.2: $p = 2$; $8 \le c \le 11$; ${\bold Z}[D]$ is non-primitive in 
$H^2({\widetilde X}, {\bold Z})$; $\pi_1(X^0) = {\bold Z}/(2)$; $\text{\rm Sing} Y = 2(c-8)A_1$.

\par \noindent
No.3: $p = 2$; $c = 12$; $Supp D$ contains only one $H$;
$\pi_1(X^0) = {\bold Z}/(2)$; $\text{\rm Sing} Y = 8A_1$.

\par \noindent
No.4: $p = 2$; $c = 12$; $Supp D = Supp(H_1 + H_2)$;
$\pi_1(X^0) = ({\bold Z}/(2))^{\oplus 2}$; $Y$ is smooth.

\par \noindent
No.5: $p = 2$; $c = 13$; ${\bold Z}[D]$ is non-primitive in 
$H^2({\widetilde X}, {\bold Z})$; $\pi_1(X^0) = ({\bold Z}/(2))^{\oplus 2}$; 
$\text{\rm Sing} Y = 4A_1$.

\par \noindent
No.6: $p = 2$; $c = 14$; ${\bold Z}[D]$ is non-primitive in 
$H^2({\widetilde X}, {\bold Z})$; $\pi_1(X^0) = ({\bold Z}/(2))^{\oplus 3}$; 
$Y$ is smooth.

\par \noindent
No.7: $p = 2$; $c = 15$; ${\bold Z}[D]$ is non-primitive in 
$H^2({\widetilde X}, {\bold Z})$; $\pi_1(X^0) = (({\bold Z}/(2))^{\oplus 4}$; 
$Y$ is smooth.

\par \noindent
No.8: $p = 2$; $c = 16$; ${\bold Z}[D]$ is non-primitive in 
$H^2({\widetilde X}, {\bold Z})$; $\pi_1(X^0)/({\bold Z}^{\oplus 4}) = {\bold Z}/(2)$; 
$Y = {\bold C}^2$.

\par \vskip 1pc \noindent
No.9: $p = 3$; $1 \le c \le 7$; ${\bold Z}[D]$ is primitive in 
$H^2({\widetilde X}, {\bold Z})$; $\pi_1(X^0) = (1)$; $Y = X$.

\par \noindent
No.10: $p = 3$; $6 \le c \le 7$; ${\bold Z}[D]$ is non-primitive in 
$H^2({\widetilde X}, {\bold Z})$; $\pi_1(X^0) = {\bold Z}/(3)$; 
$\text{\rm Sing} Y = 3(c-6)A_2$.

\par \noindent
No.11: $p = 3$; $c = 8$; $Supp D$
contains only one $R$; $\pi_1(X^0) = {\bold Z}/(3)$; 
$\text{\rm Sing} Y = 6A_2$.

\par \noindent
No.12: $p = 3$; $c = 8$; $Supp D = Supp(R_1 + R_2)$; 
$\pi_1(X^0) = (({\bold Z}/(3))^{\oplus 2}$; 
$Y$ is smooth.

\par \noindent
No.13: $p = 3$; $c = 9$; ${\bold Z}[D]$ is non-primitive in 
$H^2({\widetilde X}, {\bold Z})$; $\pi_1(X^0)/({\bold Z}^{\oplus 4}) = {\bold Z}/(3)$; 
$Y = {\bold C}^2$.

\par \vskip 1pc \noindent
No.14: $p = 5$; $1 \le c \le 4$; ${\bold Z}[D]$ is primitive in 
$H^2({\widetilde X}, {\bold Z})$; $\pi_1(X^0) = (1)$; $Y = X$.

\par \noindent
No.15: $p = 5$; $c = 4$; ${\bold Z}[D]$ is non-primitive in 
$H^2({\widetilde X}, {\bold Z})$; $\pi_1(X^0) = {\bold Z}/(5)$; 
$Y$ is smooth.

\par \vskip 1pc \noindent
No.16: $p = 7$; $1 \le c \le 3$; ${\bold Z}[D]$ is primitive in 
$H^2({\widetilde X}, {\bold Z})$; $\pi_1(X^0) = (1)$; $Y = X$.

\par \noindent
No.17: $p = 7$; $c = 3$; ${\bold Z}[D]$ is non-primitive in 
$H^2({\widetilde X}, {\bold Z})$; $\pi_1(X^0) = {\bold Z}/(7)$; 
$Y$ is smooth.

\par \vskip 1pc \noindent
No.18: $p > 7$; $c \ge 1$; ${\bold Z}[D]$ is primitive in 
$H^2({\widetilde X}, {\bold Z})$; $\pi_1(X^0) = (1)$; $Y = X$.

\par \vskip 2pc \noindent
In Table 2 below, ${\widetilde X}$ is the $K3$ cover of the Enriques
surface $W$. $W^0$ is $W$ minus a configuration $G$ of ${\bold P}^1$'s of
Dynkin type $cA_{p-1}$. $D$ is the inverse on ${\widetilde X}$ of $G$;
so $D$ has Dynkin type $2cA_{p-1}$.
$H$ or $H_i$ (resp. $K$ or $K_i$) is a 2-divisible configuration
of ${\bold P}^1$'s of Dynkin type $8A_1$ (resp. $4A_1$) on ${\widetilde X}$
(resp. on $W$).
$R$ or $R_i$ (resp. $T$) is a 3-divisible configuration of ${\bold P}^1$'s of Dynkin 
type $6A_2$ (resp. $3A_2$) on ${\widetilde X}$ (resp. on $W$).
$A_1$ is a ${\bold P}^1$ on $W$. 
$S_3$ is the symmetric group on 3 letters.
$D_{10}$ is the dihedral group of order 10.
We do not know if No.14 or No.20 in Table 2
is realizable.

\par \vskip 1pc \noindent
{\bf Table 2}

\par \vskip 1pc \noindent
No.1: $p = 2$; $1 \le c \le 3$; ${\bold Z}[G]$ is primitive in 
$H^2(W, {\bold Z})$; ${\bold Z}[D]$ is primitive in 
$H^2({\widetilde X}, {\bold Z})$; $\pi_1(W^0) = {\bold Z}/(2)$.

\par \noindent
No.2: $p = 2$; $c = 4$; ${\bold Z}[G]$ is primitive in 
$H^2(W, {\bold Z})$; ${\bold Z}[D]$ is primitive in 
$H^2({\widetilde X}, {\bold Z})$; $\pi_1(W^0) = {\bold Z}/(2)$.

\par \noindent
No.3: $p = 2$; $c = 4$; ${\bold Z}[G]$ is non-primitive in 
$H^2(W, {\bold Z})$; ${\bold Z}[D]$ is non-primitive in 
$H^2({\widetilde X}, {\bold Z})$; $\pi_1(W^0) = ({\bold Z}/(2))^{\oplus 2}$.

\par \noindent
No.4: $p = 2$; $c = 4$; ${\bold Z}[G]$ is primitive in 
$H^2(W, {\bold Z})$; ${\bold Z}[D]$ is non-primitive in 
$H^2({\widetilde X}, {\bold Z})$; $\pi_1(W^0) = {\bold Z}/(4)$.

\par \noindent
No.5: $p = 2$; $c = 5$; ${\bold Z}[G]$ is primitive in 
$H^2(W, {\bold Z})$; ${\bold Z}[D]$ is primitive in 
$H^2({\widetilde X}, {\bold Z})$; $\pi_1(W^0) = {\bold Z}/(2)$.

\par \noindent
No.6: $p = 2$; $c = 5$; ${\bold Z}[G]$ is non-primitive in 
$H^2(W, {\bold Z})$; ${\bold Z}[D]$ is non-primitive in 
$H^2({\widetilde X}, {\bold Z})$; $\pi_1(W^0) = ({\bold Z}/(2))^{\oplus 2}$.

\par \noindent
No.7: $p = 2$; $c = 5$; ${\bold Z}[G]$ is primitive in 
$H^2(W, {\bold Z})$; ${\bold Z}[D]$ is non-primitive in 
$H^2({\widetilde X}, {\bold Z})$; $\pi_1(W^0) = {\bold Z}/(4)$.

\par \noindent
No.8: $p = 2$; $c = 6$; ${\bold Z}[G]$ is primitive in 
$H^2(W, {\bold Z})$; $\pi_1(W^0) = {\bold Z}/(4)$.

\par \noindent
No.9: $p = 2$; $c = 6$; $Supp G$ contains only one $K$;
$Supp D$ contains only one $H$; $\pi_1(W^0) = ({\bold Z}/(2))^{\oplus 2}$.

\par \noindent
No.10: $p = 2$; $c = 6$; $Supp G$ contains only one $K$;
$Supp D = Supp(H_1 + H_2)$; $\pi_1(W^0) = {\bold Z}/(4) \times {\bold Z}/(2)$.

\par \noindent
No.11: $p = 2$; $c = 6$; $Supp G = Supp(K_1 + K_2)$; 
$Supp D = Supp(H_1 + H_2)$; $\pi_1(W^0) = ({\bold Z}/(2))^{\oplus 3}$.

\par \noindent
No.12: $p = 2$; $c = 7$; $Supp G$ contains only one $K$;
$Supp D = Supp(H_1 + H_2 + H_3)$; $\pi_1(W^0) = 
({\bold Z}/(4) \times {\bold Z}/(2)) \rtimes {\bold Z}/(2)$.

\par \noindent
No.13: $p = 2$; $c = 7$; $Supp G = Supp(K_1 + K_2 + A_1)$;
$Supp D = Supp(H_1 + H_2 + H_3)$; $\pi_1(W^0) = 
{\bold Z}/(4) \times ({\bold Z}/(2))^{\oplus 2}$.

\par \noindent
No.14: $p = 2$; $c = 7$; $Supp G = Supp(K_1 + K_2 + K_3)$;
$Supp D = Supp(H_1 + H_2 + H_3)$; $\pi_1(W^0) = 
({\bold Z}/(2))^{\oplus 4}$.

\par \noindent
No.15: $p = 2$; $c = 8$; ${\bold Z}[G]$ is non-primitive in
$H^2(W, {\bold Z})$; ${\bold Z}[D]$ is non-primitive in
$H^2({\widetilde X}, {\bold Z})$; $\pi_1(W^0) = 
({\bold Z}^{\oplus 2} \rtimes {\bold Z}/(2)) \rtimes {\bold Z}/(2)$.

\par \vskip 1pc \noindent
No.16: $p = 3$; $1 \le c \le 2$; ${\bold Z}[G]$ is primitive in
$H^2(W, {\bold Z})$; ${\bold Z}[D]$ is primitive in
$H^2({\widetilde X}, {\bold Z})$; $\pi_1(W^0) = {\bold Z}/(2)$.

\par \noindent
No.17: $p = 3$; $c = 3$; ${\bold Z}[G]$ is primitive in
$H^2(W, {\bold Z})$; ${\bold Z}[D]$ is primitive in
$H^2({\widetilde X}, {\bold Z})$; $\pi_1(W^0) = {\bold Z}/(2)$.

\par \noindent
No.18: $p = 3$; $c = 3$; ${\bold Z}[G]$ is non-primitive in
$H^2(W, {\bold Z})$; ${\bold Z}[D]$ is non-primitive in
$H^2({\widetilde X}, {\bold Z})$; $\pi_1(W^0) = {\bold Z}/(6)$.

\par \noindent
No.19: $p = 3$; $c = 3$; ${\bold Z}[G]$ is primitive in
$H^2(W, {\bold Z})$; ${\bold Z}[D]$ is non-primitive in
$H^2({\widetilde X}, {\bold Z})$; $\pi_1(W^0) = S_3$.

\par \noindent
No.20: $p = 3$; $c = 4$; $Supp G$ contains only one $T$;
$Supp D$ contains only one $R$; $\pi_1(W^0) = {\bold Z}/(6)$.

\par \noindent
No.21: $p = 3$; $c = 4$; $Supp G$ contains only one $T$;
$Supp D = Supp(R_1 + R_2)$; $\pi_1(W^0) = S_3 \times {\bold Z}/(3)$.

\par \vskip 1pc \noindent
No.22: $p = 5$; $c = 1$; ${\bold Z}[G]$ is primitive in
$H^2(W, {\bold Z})$; ${\bold Z}[D]$ is primitive in
$H^2({\widetilde X}, {\bold Z})$; $\pi_1(W^0) = {\bold Z}/(2)$.

\par \noindent
No.23: $p = 5$; $c = 2$; ${\bold Z}[G]$ is primitive in
$H^2(W, {\bold Z})$; ${\bold Z}[D]$ is primitive in
$H^2({\widetilde X}, {\bold Z})$; $\pi_1(W^0) = {\bold Z}/(2)$.

\par \noindent
No.24: $p = 5$; $c = 2$; ${\bold Z}[G]$ is non-primitive in
$H^2(W, {\bold Z})$; ${\bold Z}[D]$ is non-primitive in
$H^2({\widetilde X}, {\bold Z})$; $\pi_1(W^0) = {\bold Z}/(10)$.

\par \noindent
No.25: $p = 5$; $c = 2$; ${\bold Z}[G]$ is primitive in
$H^2(W, {\bold Z})$; ${\bold Z}[D]$ is non-primitive in
$H^2({\widetilde X}, {\bold Z})$; $\pi_1(W^0) = D_{10}$.

\par \vskip 1pc \noindent
No.26: $p = 7$; $c = 1$; ${\bold Z}[G]$ is primitive in
$H^2(W, {\bold Z})$; ${\bold Z}[D]$ is primitive in
$H^2({\widetilde X}, {\bold Z})$; $\pi_1(W^0) = {\bold Z}/(2)$.

\enddocument